\documentclass[article]{siamart250211}
\usepackage[caption=false]{subfig}
\usepackage{nicematrix}
\title{Fast solver for the Reynolds Equation \linebreak on piecewise linear geometries}

\author{Sarah Dennis\thanks{Department of Mathematics, Brandeis University, Waltham MA (\email{sarahdennis@brandeis.edu})} \and Thomas G. Fai\thanks{Department of Mathematics and Volen Center for Complex Systems, Brandeis University, Waltham MA (\email{tfai@brandeis.edu})}}

\begin{document}
\maketitle
\begin{abstract}
The Reynolds equation is derived from the incompressible Navier Stokes equations under the lubrication assumptions of a long and thin domain geometry and a small scaled Reynolds number. The Reynolds equation is an elliptic differential equation and a dramatic simplification from the governing equations. When the fluid domain is piecewise linear, the Reynolds equation has an exact solution that we formulate by coupling the exact solutions of each piecewise component. We consider a formulation specifically for piecewise constant heights, and a more general formulation for piecewise linear heights; in both cases the linear system is inverted using the Schur complement. These methods can also be applied in the case of non-linear heights by approximating the height as piecewise constant or piecewise linear, in which case the methods achieve second-order accuracy. We assess the time complexity of the two methods, and determine that the method for piecewise linear heights is linear time for the number of piecewise components. As an application of these methods, we explore the limits of validity for lubrication theory by comparing the solutions of the Reynolds and the Stokes equations for a variety of linear and non-linear textured slider geometries.
\end{abstract}
\section{Introduction}
The Reynolds equation of lubrication theory is an elliptic differential equation derived from the incompressible Navier-Stokes equations under the lubrication assumptions of a long, thin fluid domain and a small scaled Reynolds number. In contrast, the Stokes equation is derived from the incompressible Navier-Stokes equations under the sole assumption of zero Reynolds number. That is, Stokes flow has relaxed restrictions on the ratio of length scales compared to lubrication theory, providing a standpoint from which we observe the sensitivity of the Reynolds equation to the thin film condition. 

Because lubrication theory is based in the assumption of a long and thin fluid domain, solutions in this regime are highly sensitive to variations in the height of the fluid film. In the limit of zero Reynolds number flows, discrepancies between lubrication theory and Stokes flow are found to increase not only with increasing ratio of the length scales, but also with increasing magnitude of gradients in the fluid film height \cite{brown_applicability_1995,biswas_backward-facing_2004,dobrica_about_2009}. Accordingly, it may be that the fluid geometry is overall long and thin, yet surface gradients of large magnitude also impact the extent to which the thin film assumption is satisfied, thus contributing to error in the lubrication theory solution compared with Stokes flow. When the lubrication assumptions break down, the total pressure drop modeled by lubrication theory is found to be an underestimate compared with the total pressure drop modeled by Stokes flow \cite{dobrica_reynolds_2005,dobrica_about_2009}. Furthermore, lubrication theory does not accurately capture flow separation and flow recirculation resulting from large surface gradients or corners in the fluid geometry \cite{armaly_experimental_1983,biswas_backward-facing_2004,shyu_numerical_2018}. Here we consider a variety of textured sliders, including several which feature large surface gradients and sharp corners, so as to further analyze these differences between the solutions of lubrication theory and Stokes flow.

When the fluid domain has a linear or constant height profile, the Reynolds equation has an exact solution. Approximating the height as piecewise constant or piecewise linear, we couple the exact solutions to the Reynolds equation on each piecewise component using the conditions of constant flux and continuous pressure. The resulting linear system of equations can be efficiently solved using the Schur complement \cite{gallier_geometric_2011}, leading to an exact solution for the Reynolds equation in the case of piecewise linear heights, and a second-order accurate solution for arbitrary non-linear heights. Compared with a standard finite difference approach, the piecewise analytic methods are robust, retaining second-order accuracy even for discontinuous height profiles. As an application of the methods for the Reynolds equation, we explore the range of validity for lubrication theory by comparing to Stokes flow. 

We present two methods of solution for the Reynolds equation: one approach considers a piecewise constant approximation of the fluid domain, and the other considers a more general piecewise linear approximation. A similar approach of piecewise exact solutions to the Reynolds equation is considered in \cite{rahmani_analytical_2010}; here, we also assess the performance of the piecewise linear (PWL) and piecewise constant (PWC) methods in approximating the solution for non-linear height functions. We compare the performance of each method for the Reynolds equation in terms of time complexity, including a comparison to a standard finite difference (FD) method; we determine that the PWL method performs fastest, running in linear time for the number of piecewise components. Although the PWL method is more general and has better time complexity, the PWC method has a simpler formulation through which we introduce the process of solving the Reynolds equation in piecewise components.  

Finally, to demonstrate the PWL and PWC methods and to examine the range of validity for lubrication theory, we consider four examples of textured sliders: a backward facing step, a wedge slider, a logistic step, and a sinusoidal slider. We contrast the solutions from the Reynolds equation and the Stokes equation, confirming that large surface gradients cause the lubrication assumptions to break down and correspond to significant discrepancies between lubrication theory and Stokes flows.

Source code is available at \href{https://github.com/sarah-dennis/piecewise-reynolds}{https://github.com/sarah-dennis/piecewise-reynolds}.
\section{Lubrication Theory}
\begin{figure}
    \centering
    \includegraphics[width=0.95\textwidth]{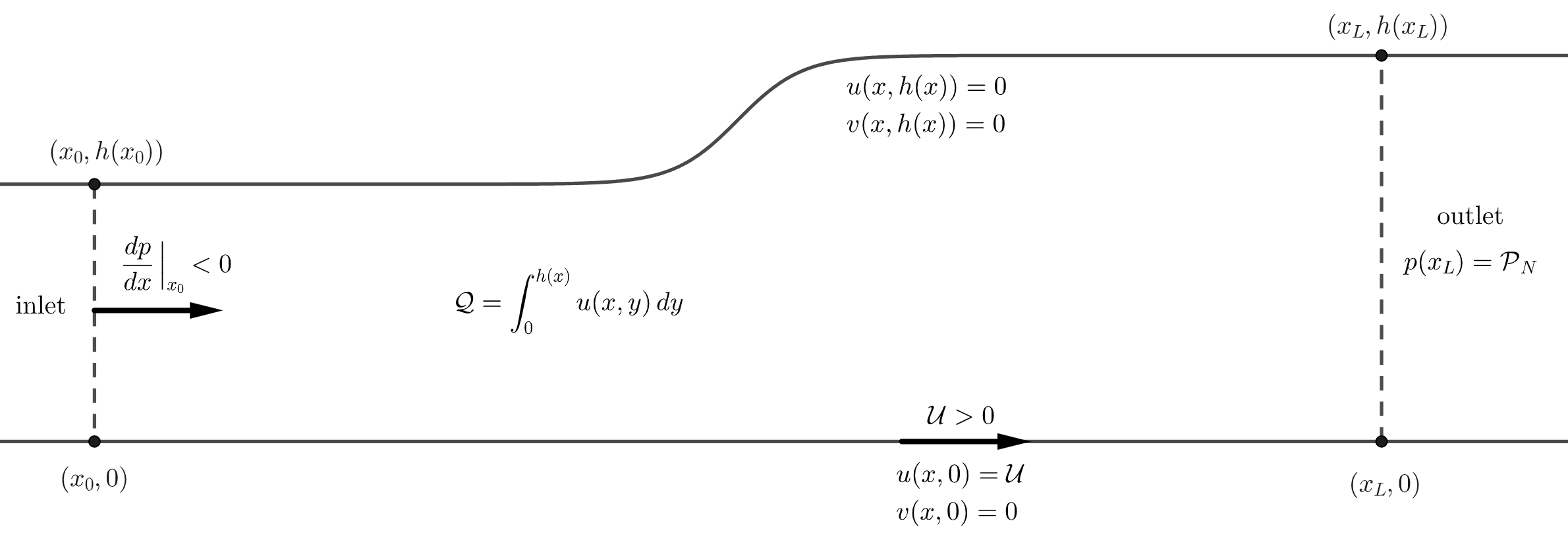}
    \caption{\color{black}Schematic of the two dimensional fluid domain with inlet $x=x_0$, outlet $x=x_L$ and height $y=h(x)$. The constant flux $\mathcal{Q}$, boundary velocity $\mathcal{U}$ and outlet pressure $\mathcal{P}_N$ are prescribed.}
    \label{schematic_bcs}
\end{figure}

The Reynolds equation is derived from the Navier-Stokes equations under the lubrication assumptions \cite{leal_advanced_2007}. For the two dimensional fluid domain $(x,y) \in [x_0,x_L]\times[0,h(x)]$, {\color{black}shown in \cref{schematic_bcs}}, denote the characteristic length scales $L_x = x_L - x_0$ and $L_y = \max h(x) > 0$, and the length scale ratio $\varepsilon = L_y/L_x$. Given a prescribed constant flux $\mathcal{Q}\ne 0$, the characteristic fluid velocities are given by $U_* = \mathcal{Q}/L_y$ and $V_* = \mathcal{Q}/L_x$. The characteristic time is $T_* = L_x/U_*$ and the characteristic pressure is $P_* = \eta U_*L_x/L_y^2$.
The Reynolds number is given by $\text{Re}=\rho U_*L_x/\eta$, where $\eta$ is the constant dynamic viscosity ($\eta = 1$) and $\rho$ is the constant density ($\rho = 1$). 
The lubrication assumptions $\varepsilon \ll 1$ and $\varepsilon^2 \text{Re}\ll 1$ characterize a long and thin fluid with a small scaled Reynolds number; these assumptions yield an approximation to the Navier-Stokes equations \cref{n-s_x,n-s_y}, resulting in the momentum equations,
\begin{align}\label{gap_x} 
 \frac{\partial p}{\partial x}&=\eta\frac{\partial^2 u}{\partial y^2}, \\ \label{gap_y}
 \frac{\partial p}{\partial y}&=0.
\end{align}
Together with incompressibility, \begin{equation}\label{incompressibility}
 \frac{\partial u}{\partial x} + \frac{\partial v}{\partial y} = 0,\end{equation} \cref{gap_x,gap_y} constitute the governing equations for lubrication theory.

{\color{black} A schematic for the two dimensional fluid domain is shown in \cref{schematic_bcs}.} We assume the no-slip boundary condition at the fluid-surface interfaces $y=0$ and $y=h(x)$. Without loss of generality, the velocity boundary conditions are,
\begin{align}\label{reyn_bc_u}
&& u (x,0) = \mathcal{U}, && u(x,h(x)) = 0, \\\label{reyn_bc_v}
&& v (x,0) = 0, && v(x,h(x)) = 0.
\end{align}

The velocity $u(x,y)$ is determined through integration of \cref{gap_x} and applying the boundary conditions \cref{reyn_bc_u},
\begin{equation}
 \label{reyn_u} u(x,y)=\frac{1}{2\eta}\frac{d p}{d x} \Big(y^2-h(x)y\Big) +\frac{\mathcal{U}}{h(x)} \Big(h(x)-y\Big).\end{equation} The velocity $v(x,y)$ is then determined from incompressibility \cref{incompressibility} and applying the boundary conditions \cref{reyn_bc_v},
\begin{equation}\label{reyn_v}
v(x,y) =\frac{-1}{6\eta}\frac{d^2 p}{d x^2}y^3 + \frac{1}{2}\Bigg(\frac{1}{2\eta}\bigg(\frac{d^2 p}{d x^2}h(x) + \frac{d p}{d x}\frac{dh}{dx}\bigg) -\frac{\mathcal{U}}{[h(x)]^2}\frac{dh}{dx}\Bigg)y^2.\end{equation} 
From incompressibility, the condition of constant flux $\mathcal{Q}$,
\begin{equation}\label{flux}\mathcal{Q} = \int_{0}^{h(x)} u(x,y) \, dy =\frac{-1}{12\eta}\Bigg(\Big[h(x)\Big]^3\frac{dp}{dx} -6\eta\mathcal{U}h(x)\Bigg),\end{equation} and the boundary condition $v(x,h)=0$, are satisfied exactly when $p(x)$ satisfies the Reynolds equation,
\begin{equation}\label{reynolds}
 \frac{d}{dx}\Bigg[\Big[h(x)\Big]^3\frac{dp}{dx}\Bigg] = 6\eta\,\mathcal{U}\frac{dh}{dx}.
\end{equation}

For the pressure boundary conditions, we consider a mixed Dirichlet-Neumann boundary condition, prescribing the flux $\mathcal{Q}$ and outlet pressure $\mathcal{P}_N$,
\begin{align}\label{reyn_bc_p_mixed}&& \frac{dp}{dx}\bigg|_{x_0} = \frac{-12\eta \mathcal{Q}}{\big[h(x_0)\big]^3} +\frac{6\eta\mathcal{U}}{\big[h(x_0)\big]^2}, &&p (x_L) =\mathcal{P}_N. &&\end{align} {\color{black}A negative pressure gradient is prescribed at the inlet, $\frac{dp}{dx}\big|_{x_0} < 0$, so that the bulk flow travels in the direction of the outlet $x=x_L$.}
To translate \cref{reyn_bc_p_mixed} to a Dirichlet boundary condition, the prescribed flux $\mathcal{Q}$ is exchanged for a prescribed inlet pressure $p(x_0)=\mathcal{P}_0$, corresponding to a prescribed pressure drop $\Delta p = p(x_0) - p(x_L)$ for the whole domain,
\begin{equation}\label{dP} \Delta p= -12\eta \mathcal{Q}\int_{x_L}^{x_0}\big[h(x)\big]^{-3}dx + 6\eta\mathcal{U}\int_{x_L}^{x_0} \big[h(x)\big]^{-2}dx.
\end{equation}
In the case that height $h(x)$ is linear, the integrals involving $h(x)$ in \cref{dP} can be evaluated exactly. Evaluating these integrals is also the key to obtaining an exact solution for the Reynolds equation.

\section{The piecewise constant solution} The following method (PWC) for the Reynolds equation considers a piecewise constant approximation of the height function, and utilizes the Schur complement factorization technique to solve the linear system. 
In the case that $h(x)>0$ is constant, the Reynolds equation \cref{reynolds} reduces to $\frac{d^2p}{dx^2} = 0$ and the solution $p(x)$ is linear. 

Consider a piecewise constant height $h(x)$ with $N$ components. Let $\{x_k\}_0^N$ demarcate the piecewise constant regions of $h(x)$; denote the widths $\Delta x_k=x_{k+1}-x_k$. Let $\{h_k\}_0^{N-1}$ denote the constant value of $h(x)$ on $(x_k, x_{k+1})$, and let $\{\frac{\Delta p}{\Delta x}\big|_k\}_0^{N-1}$ denote the constant pressure gradient on $(x_k, x_{k+1})$, 
\begin{equation} \label{pwc_dp} \frac{\Delta p}{\Delta x}\Big|_k = \frac{p_{k+1}-p_{k}}{x_{k+1}-x_k},\end{equation} 
where $\{p_k\}_0^N$ are the pressure endpoints,
$p_k = p(x_k)$. The flux $\mathcal{Q}$ expressed in \cref{flux} relates each of the $N$ constant pressure gradients. Eliminating the constant flux gives the relationship,
\begin{equation}\label{pwc_dq}
h_k^3\frac{\Delta p}{\Delta x}\Big|_k - h_{j}^3\frac{\Delta p}{\Delta x}\Big|_{j}  + 6\eta \mathcal{U} (h_{j}- h_k)=0,
\end{equation}for any $j,k \in [0,N-1]$. 

 The expressions \cref{pwc_dp} for $k\in [0,N-1]$ and \cref{pwc_dq} for $j=k+1$, $k\in [0,N-2]$ constitute a size $2N-1$ linear system $\mathrm{M}{\bf x}={\bf b}$ solving for, \begin{equation}\label{pwc_x}{\bf x} = \bigg[\frac{\Delta p}{\Delta x}\Big|_0, \frac{\Delta p}{\Delta x}\Big|_1, \dots \frac{\Delta p}{\Delta x}\Big|_{N-1}, p_1, p_2, \dots p_{N-1}\bigg]^T.\end{equation}
The matrix $\mathrm{M}$ is formulated as a block matrix,
\begin{equation}\label{pwc_linsys}
\mathrm{M} = \begin{bNiceArray}{ccccc|cccc}[margin,columns-width=.25cm]
  \Block{5-5}<\Large>{\mathcal{I}} & & & & & 0 & 0 & \cdots & 0 \\
  & & & & &                         \frac{1}{\Delta x_1} & \frac{-1}{\Delta x_1} & \ddots & \vdots \\
  & & & & &                          0 &\ddots & \ddots & 0 \\
  & & & & &                          \vdots & \ddots & \frac{1}{\Delta x_{N-2}} & \frac{-1}{\Delta x_{N-2}} \\[.5em]
  & & & & &                          0 & \cdots & 0 & \frac{1}{\Delta x_{N-1}}
  \\[.5em]\hline\\[-.5em]
 -h_0^3 & h_1^3 & 0 & \cdots & 0 & \Block{4-4}<\Large>{0}& & & \\
 0 & -h_1^3 & h_2^3 & \ddots & \vdots & & & & \\
 \vdots & \ddots & \ddots & \ddots & 0 & & & &\\
 0 & \cdots & 0 & -h_{N-2}^3 & h_{N-1}^3  & & & &\\
\end{bNiceArray},\end{equation}
where $\mathcal{I}$ denotes an identity block of size $N$. The right hand side is given by, 
\begin{multline}\label{pwc_rhs}{\bf b} = \bigg[-12\eta\Big(\mathcal{Q}h_0^{-3}-\tfrac{1}{2}\mathcal{U}h_0^{-2}\Big), 0, \dots, 0,  \tfrac{\mathcal{P}_N}{\Delta x_{N-1}},\\ 6\eta\mathcal{U}(h_1-h_0), 6\eta\mathcal{U} (h_2-h_1), \dots, 6\eta\mathcal{U}(h_{N-1}-h_{N-2})\bigg]^T.\end{multline}

{\color{black}To evaluate the linear system $\mathrm{M}{\bf x}={\bf b}$ for \cref{pwc_x,pwc_linsys,pwc_rhs}, we utilize the Schur complement.} 
The Schur complement is a linear algebra technique utilized in the factorization and inversion of block matrices. The theory is presented in \cite{gallier_geometric_2011} and summarized here.
Given a block matrix,
\begin{equation}\label{block_M} \mathrm{M} = \begin{bmatrix}
    A & B \\ C & D
\end{bmatrix},\end{equation}
where $A$ is invertible, the Schur complement of $\mathrm{M}$ is given by,
\begin{equation}\label{schur} K = D - CA^{-1}B.\end{equation}
Through Gaussian elimination, the Schur complement $K$ gives a block LDU decomposition for $\mathrm{M}$,
\begin{equation}\mathrm{M} = \begin{bmatrix}
    \mathcal{I} & 0 \\ CA^{-1}  & \mathcal{I}
\end{bmatrix}
 \begin{bmatrix}
    A& 0 \\0  & K
\end{bmatrix}
 \begin{bmatrix}
    \mathcal{I} & A^{-1}B \\0 & \mathcal{I}
\end{bmatrix},
\end{equation} and reduces the problem of finding $\mathrm{M}^{-1}$ to finding $K^{-1}$ and $A^{-1}$. Following the process of forward and backward substitution from the block LDU decomposition for $\mathrm{M}$, the general form of $\mathrm{M}^{-1}$ is derived,
\begin{equation} \mathrm{M}^{-1} = \begin{bmatrix}
A^{-1} + A^{-1}BK^{-1}CA^{-1} & -A^{-1}BK^{-1} \\ 
-K^{-1}CA^{-1} & K^{-1} 
\end{bmatrix}.\end{equation} 
{\color{black}Hence, the particular matrix $\mathrm{M}$ in \cref{pwc_linsys} has inverse,}
\begin{equation}\label{M_inv} \mathrm{M}^{-1} = \begin{bmatrix}
\mathcal{I} + BK^{-1}C & -BK^{-1} \\ 
-K^{-1}C & K^{-1}
\end{bmatrix},\end{equation} where $K$ is the Schur complement of size $N-1$.

In the implementation, one need not compute all entries to the matrix products appearing in the blocks of $\mathrm{M}^{-1}$, although it is necessary to compute all entries to the dense matrix $K^{-1}$. Observe that at most two entries in the upper block of ${\bf b}$ are non-zero, corresponding to the inlet and outlet boundary conditions for pressure. Then on evaluating $\mathrm{M}^{-1}{\bf b} = {\bf x}$ to solve for $\{p_k\}_{k=1}^{N-1}$, it is only necessary to compute the first and last columns of $-K^{-1}C$. Once the pressure endpoints $\{p_k\}_{k=1}^{N-1}$ are determined from the upper matrix blocks, the pressure gradients $\{\frac{\Delta p}{\Delta x}\}_{k=0}^{N-1}$ are easily obtained. In this sense, the solution hinges almost entirely on inverting the Schur complement.

In this case, the Schur complement for $\mathrm{M}$ is simply $K = -CB$, which is a symmetric tri-diagonal matrix, with off-diagonal elements $K_{i,i+1} = K_{i+1,i} =\frac{h_{i+1}^3}{\Delta x_{i+1}}$ and diagonal elements $K_{i,i} = -\Big(\frac{h_{i}^3}{\Delta x_i}+\frac{h_{i+1}^3}{\Delta x_{i+1}}\Big)$. 
An efficient and numerically stable algorithm for the inversion of tri-diagonal matrices is presented in \cite{jain_numerically_2007}, and here we outline the reduced algorithm for the symmetric case.

Let $K_{i,i} = \alpha_i$ for $i\in [0, N-2]$ and $K_{i,i+1} = K_{i+1,i} = \beta_i$ for $i\in [0,N-3]$ denote the diagonal and off-diagonal elements of the Schur complement $K$. Define the recursive sequence $\{S_k\}_{N-3}^0$ as,
\begin{equation}
    \begin{cases}
    S_{N-3} = -\beta_{N-3}/\alpha_{N-2} \\
    S_{k} = -\beta_k/(\alpha_{k+1} + S_{k+1}\beta_{k+1}) & k\in [0, N-4]
\end{cases}.
\end{equation}
The diagonal elements of $K^{-1}$ are then obtained recursively as,
\begin{equation}
    \begin{cases}
K^{-1}_{0,0} = 1/(\alpha_0 + \beta_0 S_0)\\
K^{-1}_{i+1,i+1} = (1 - \beta_{i}K^{-1}_{i,i}S_i)/(\alpha_{i+1} +\beta_{i+1}S_{i+1}) & i\in [0, N-3]\\
K^{-1}_{N-2,N-2} = (1-\beta_{N-3}K^{-1}_{N-3,N-3}S_{N-3})/\alpha_{N-2}
\end{cases},
\end{equation}
and the remaining elements of the symmetric $K^{-1}$ are given by,
\begin{equation}
    \begin{cases}
    K^{-1}_{i,j} = K^{-1}_{i,i}\prod_{k=i}^{j-1}S_k & i < j\\
    K^{-1}_{i,j} = K^{-1}_{j,i}& i > j 
\end{cases}.
\end{equation}

Through this approach, we gain element-wise access to $K^{-1}$ through two recursive sequences $\{S_k\}_{N-3}^0$ and $\{K^{-1}_{i,i}\}_{0}^{N-2}$. In the implementation, we also precompute the partial product sequence $\{T_k\}_{k=0}^{N-3}=\{\Pi_{i=0}^k S_k\}_{k=0}^{N-3}$ so that non-diagonal elements of $K^{-1}_{i<j}$ are evaluated as $K^{-1}_{i,j} = K^{-1}_{i,i} \big(T_{j-1} / T_{i-1}\big)$. Once $K^{-1}$ is determined, the solution $M^{-1}{\bf b}={\bf x}$ can be evaluated. Although $K^{-1}$ can be described in $\mathcal{O}(N)$ time with the two recursive sequences, evaluating $M^{-1}{\bf b}={\bf x}$ is done in $\mathcal{O}(N^2)$ time. 

\section{The piecewise linear solution} The following method (PWL) for the \linebreak Reynolds equation considers a piecewise linear approximation to the height, and formulates a new linear system that can also be solved using the Schur complement. The PWL method is more general and more efficient than the PWC method. Where the PWC method uses the constant pressure gradients $\frac{\Delta p}{\Delta x}\big|_k$ as variables in the linear system, the PWL method uses a single variable $\mathcal{C}_Q = -12\eta \mathcal{Q}$ corresponding to the flux, unifying the (now not necessarily constant) pressure gradients $\frac{dp}{dx}\big|_k$ on each component. 

Suppose $h(x)$ defined on $[x_0,x_L]$ is piecewise linear with $N$ components. Let $\{x_k\}_0^N$ denote the critical points of $h(x)$ demarcating each piecewise linear region. Let $\{h_{k^\pm}\}_{k=0}^N$ denote the endpoints of $h(x)$ on each interval, $h_{k^\pm} = \lim_{x\to x_k^{\pm}}h(x)$. And let $\{\frac{\Delta h}{\Delta x}\big|_k\}_{k=0}^{N-1}$ denote the constant gradients of $h(x)$ for each interval, \begin{align}
    \frac{\Delta h}{\Delta x}\bigg|_k = \frac{h_{(k+1)^-}-h_{k^+}}{x_{k+1}-x_k} = \lim_{x\to x_k^+}\frac{dh}{dx}.
\end{align}
The solution to the Reynolds equation for a piecewise linear height is derived from coupling solutions of the form \cref{pwl_p} on each sub-interval. 

For a single interval $[x_k,x_{k+1}]$ on which $h(x)$ is piecewise linear, $p(x)$ satisfies,
\begin{equation}\label{pwl_p}
    p(x) = \begin{cases}
               -\Bigg(\tfrac{1}{2}\mathcal{C}_Q \Big[h(x)\Big]^{-2} + 6 \eta \mathcal{U} \Big[h(x)\Big]^{-1}\Bigg)\Big[\frac{\Delta h}{\Delta x}\Big|_k\,\Big]^{-1}  + \mathcal{C}_{P_k} &\frac{\Delta h}{ \Delta x}\Big|_k \ne 0\\
    \Big(\mathcal{C}_Q \Big[h(x)\Big]^{-3}+ 6 \eta \mathcal{U} \Big[h(x)\Big]^{-2}\Big)(x-x_k) + \mathcal{C}_{P_k}  & \frac{\Delta h}{\Delta x}\Big|_k = 0 \end{cases}.\end{equation}
The constant $C_Q$ arises in the first integration of the Reynolds equation and is directly proportional to the flux; comparing with \cref{flux} gives $\mathcal{C}_Q =-12\eta \mathcal{Q}$. The constant $C_{P_k}$ arises in the second integration of the Reynolds equation and relates to the fixed endpoint pressure. For $h(x)$ with $N$ piecewise linear components, we couple solutions of the form \cref{pwl_p}. The fixed outlet pressure $ p(x_L)=\mathcal{P}_N$ determines $\mathcal{C}_{P_{N-1}}$,
\begin{equation}\label{boundary_cp}\mathcal{C}_{P_{N-1}} = \begin{cases}  
\mathcal{P}_N  +\mathcal{C}_Q\tfrac{1}{2}\frac{\Delta h}{\Delta x}\Big|_{N-1}^{-1}   h_{N^-}^{-2} +6 \eta \mathcal{U}\frac{\Delta h}{\Delta x} \Big|_{N-1}^{-1} h_{N^-}^{-1} &\frac{\Delta h}{ \Delta x}\Big|_{N-1}\ne 0\\
\mathcal{P}_N  -\mathcal{C}_Q  h_{N^-}^{-3}\Delta x|_{N-1} -6 \eta \mathcal{U}  h_{N^-}^{-2}\Delta x|_{N-1} & \frac{\Delta h}{\Delta x}\Big|_{N-1} = 0
\end{cases}.\end{equation} 
To solve for $\{\mathcal{C}_{P_k}\}_{k=0}^{N-2}$, we assume a continuous pressure $p(x)$ and set the left and right limits equal at each $x_k$,
\begin{align}
    \lim_{x\to x_k^-}p(x) &=
    \begin{cases} 
    -\frac{\Delta h}{\Delta x}\Big|_{k-1}^{-1} \Big(\frac{1}{2}\mathcal{C}_Q h_{k^-}^{-2} + 6 \eta \mathcal{U} h_{k^-}^{-1}\Big) + \mathcal{C}_{P_{k-1}} 
    &\frac{\Delta h}{\Delta x}\big|_{k-1} \ne 0 \\
     \Delta x|_{k-1}  \Big(\mathcal{C}_Q h_{k^-}^{-3}+ 6 \eta \mathcal{U} h_{k^-}^{-2}\Big) + \mathcal{C}_{P_{k-1}} 
     & \frac{\Delta h}{\Delta x}\big|_{k-1}= 0
    \end{cases},
\\
    \lim_{x\to x_k^+}p(x) &= 
    \begin{cases}
    -\frac{\Delta h}{\Delta x}\Big|_{k}^{-1} \Big(\frac{1}{2}\mathcal{C}_Q h_{k^+}^{-2} + 6 \eta \mathcal{U} h_{k^+}^{-1}\Big) + \mathcal{C}_{P_{k}} 
    & \frac{\Delta h}{\Delta x}\big|_{k} \ne 0 \\
  \mathcal{C}_{P_{k}} 
  & \frac{\Delta h}{\Delta x}\big|_{k}= 0
\end{cases}.
\end{align}
Hence for each $0< k < N$, the integration constants $C_{P_k}$, $C_{P_{k-1}}$ and $C_Q$ satisfy,
\begin{multline} \label{cp cq cases}
    \begin{cases}
    \mathcal{C}_{P_{k}}- \mathcal{C}_{P_{k-1}} -\mathcal{C}_Q\bigg(\frac{1}{2}h_{k^+}^{-2}\frac{\Delta h}{\Delta x}\Big|_{k}^{-1} - \frac{1}{2} h_{k^-}^{-2} \frac{\Delta h}{\Delta x}\Big|_{k-1}^{-1} \bigg) \\
        \hspace{2em} = 6\eta \mathcal{U}\bigg(h_{k^+}^{-1}\frac{\Delta h}{\Delta x}\Big|_{k}^{-1} -h_{k^-}^{-1} \frac{\Delta h}{\Delta x}\Big|_{k-1}^{-1}   \bigg) 
            & \frac{\Delta h}{\Delta x}\big|_{k} \ne 0,   \frac{\Delta h}{\Delta x}\big|_{k-1} \ne 0  \\
   \mathcal{C}_{P_{k}}- \mathcal{C}_{P_{k-1}} -\mathcal{C}_Q\bigg(\frac{1}{2}h_{k^+}^{-2} \frac{\Delta h}{\Delta x}\Big|_{k}^{-1}  + h_{k^-}^{-3}\Delta x\Big|_{k-1}\bigg) \\
        \hspace{2em} = 6\eta \mathcal{U}\bigg( h_{k^+}^{-1}\frac{\Delta h}{\Delta x}\Big|_{k}^{-1}  + h_{k^-}^{-2}\Delta x\Big|_{k-1}  \bigg)
            & \frac{\Delta h}{\Delta x}\big|_{k} \ne 0,  \frac{\Delta h}{\Delta x}\big|_{k-1}= 0 \\
   \mathcal{C}_{P_{k}}- \mathcal{C}_{P_{k-1}} +\mathcal{C}_Q \frac{1}{2}h_{k^-}^{-2}\frac{\Delta h}{\Delta x}\Big|_{k-1}^{-1}\\
        \hspace{2em}= -6\eta \mathcal{U}h_{k^-}^{-1}\frac{\Delta h}{\Delta x}\big|_{k-1}^{-1}   
            &  \frac{\Delta h}{\Delta x}\big|_{k} = 0,\frac{\Delta h}{\Delta x}\big|_{k-1}\ne 0\vspace{.5em}\\
    \mathcal{C}_{P_{k}}- \mathcal{C}_{P_{k-1}} -\mathcal{C}_Q  h_{k^-}^{-3} \Delta x \big|_{k-1}\\
        \hspace{2em}= 6\eta \mathcal{U} h_{k^-}^{-2}\Delta x\big|_{k-1}  
            & \frac{\Delta h}{\Delta x}\big|_{k} = 0, \frac{\Delta h}{\Delta x}\big|_{k-1} = 0
    \end{cases}.
\end{multline}

All together, $C_Q$ and $\{C_{P_k}\}_{k=0}^{N-1}$ constitute a size $N+1$ linear system of equations $\mathrm{M}{\bf x} = {\bf b}$ where ${\bf x} = [ \mathcal{C}_Q,\mathcal{C}_{P_0}, ...,  \mathcal{C}_{P_{N-1}} ]^T$. Corresponding to \cref{cp cq cases}, the matrix $\mathrm{M}$ stores the coefficients on  $C_Q$, $\mathcal{C}_{P_k}$, $\mathcal{C}_{P_{k-1}}$, and the right hand side vector ${\bf b}$ holds constant terms. 
For example, if $\frac{\Delta h}{\Delta x}\big|_{k} \ne 0$ for all $k$, 
\begin{equation}
    \mathrm{M} =     \begin{bmatrix}
    1 & 0 & & \cdots & &0  \\
    -\tfrac12 \Big(h_{1^+}^{-2}\frac{\Delta h}{\Delta x}\big|_{1}^{-1} -h_{1^-}^{-2}\frac{\Delta h}{\Delta x}\big|_{0}^{-1} \Big)&-1 &1 &  0 &\cdots&0\\
    \vdots& &\ddots &\ddots &  &  \vdots\\
     -\tfrac12 \Big(h_{k^+}^{-2}\frac{\Delta h}{\Delta x}\big|_{k}^{-1} -h_{k^-}^{-2}\frac{\Delta h}{\Delta x}\big|_{k-1}^{-1} \Big)&0&-1 &1 &  0 &\cdots\\
     \vdots& \vdots&   &\ddots   &\ddots & \\
     -\tfrac12 \Big(h_{N-1^+}^{-2}\frac{\Delta h}{\Delta x}\big|_{N-1}^{-1} -h_{N-1^-}^{-2}\frac{\Delta h}{\Delta x}\big|_{N-2}^{-1} \Big)&0&\cdots& 0 & -1 &1  \\
      \tfrac12 h_N^{-2}\frac{\Delta h}{\Delta x}\big|_{N-1}^{-1}&0&\cdots&  &0&-1
    \end{bmatrix},
\end{equation}
and,
\begin{multline} 
    {\bf b} = \bigg[
    -12\eta\mathcal{Q},
  6\eta \mathcal{U} \Big(h_{1^+}^{-1}\frac{\Delta h}{\Delta x}\Big|_1^{-1} - h_{1^-}^{-1}\frac{\Delta h}{\Delta x}\Big|_0^{-1}\Big),
    \dots,\\
    6\eta \mathcal{U} \Big(h_{k^+}^{-1}\frac{\Delta h}{\Delta x}\Big|_k^{-1}  - h_{k^-}^{-1}\frac{\Delta h}{\Delta x}\Big|_{k-1}^{-1} \Big),
     \dots,
-6\eta \mathcal{U}h_N^{-1}\frac{\Delta h}{\Delta x}\Big|_{N-1}^{-1} - \mathcal{P}_N\Bigg]^T.\end{multline}

The matrix $\mathrm{M}$ has a block form; in view of the structure \cref{block_M}, $A=1$, $B=0$, $C$ is the vector of $C_Q$ coefficients from \cref{cp cq cases}, and $D$ is a bi-diagonal matrix corresponding to the $C_{P_k}$ and $C_{P_{k-1}}$ coefficients from \cref{cp cq cases}. Since $B=0$, the Schur complement \cref{schur} is simply $K=D$, and the inverse of $\mathrm{M}$ is, 
\begin{equation}\mathrm{M}^{-1} = \begin{bmatrix}
    1 & 0\\
    -D^{-1}C & D^{-1}
\end{bmatrix}.\end{equation} 
The matrix $K^{-1}=D^{-1}$ is upper triangular with all $-1$ entries. The matrix vector product $-D^{-1}C$ is then the reversed partial sum of entries in $C$. Likewise, when evaluating $M^{-1}{\bf b}$, the lower block of ${\bf b}$ is multiplied by $D^{-1}$, resulting in the negative reversed partial sum of these entries in ${\bf b}$. In this case, solving with the Schur complement is equivalent to reducing the identity row in $\mathrm{M}$: solving $C_Q$ and moving the block $C$ to the right hand side. After solving for ${\bf x}$ comprising $C_Q$ and $\{C_{P_k}\}_{k=0}^{N-1}$, the pressure $p(x)$ on each region $\{[x_k,x_{k+1}]\}_{k=0}^{N-1}$ is obtained according to \cref{pwl_p}.

The PWL method for the Reynolds equation is more efficient than the less general PWC method. In the PWC method, we independently compute each entry to $K^{-1}$ when we evaluate $\mathrm{M^{-1}}{\bf b}={\bf x}$, leading to an $\mathcal{O}(N^2)$ time algorithm. In the PWL method, $K^{-1}$ is constant and upper triangular, so $\mathrm{M^{-1}}{\bf b}={\bf x}$ can be evaluated through computing two partial sums in an $\mathcal{O}(N)$ time algorithm. 

\section{Examples}
We now consider a variety of textured slider examples to evaluate the PWL and PWC methods for the Reynolds equation. We compare the PWC and PWL solutions with a finite difference (FD) solution for the Reynolds equation (\cref{app_reyn_fd}), and with a finite difference solution to the Stokes equation (\cref{app_stokes}). We examine both piecewise linear and non-linear height functions: the backward facing step, the wedge slider, the logistic step, and the sinusoidal slider. For the latter two examples, the PWC and PWL methods consider approximations to the non-linear height. These examples feature surface discontinuities and large surface gradients to showcase the differences in the pressure and velocity solutions from the Reynolds and the Stokes equations.

{\color{black}In each example to follow, we fix the domain length at $L_x = 16$ and consider domain heights $L_y \simeq 2$. This gives a length scale ratio $\varepsilon \simeq 1/8$, consistent with the assumptions of lubrication theory. For domains of this size, results are presented with $N = 160$, where $N=1/\Delta x = 1/\Delta y$ is the number of grid points per unit length. The 2D finite difference method for the Stokes equation utilizes $N^2L_yL_x$ total grid points. For the Reynolds equation, the 1D finite difference method utilizes $NL_x$ total grid points; likewise the PWC and PWL methods use $NL_x-1$ piecewise components in the case of arbitrary geometries. While the 1D methods may be simulated at much finer grid resolutions, the combination $N=160$ and $L_yL_x=32$ nears the memory limit of the 2D simulations given our computational resources.}

\subsection{Piecewise constant height}
\begin{figure}
    \centering
    \includegraphics[width=0.9\textwidth]{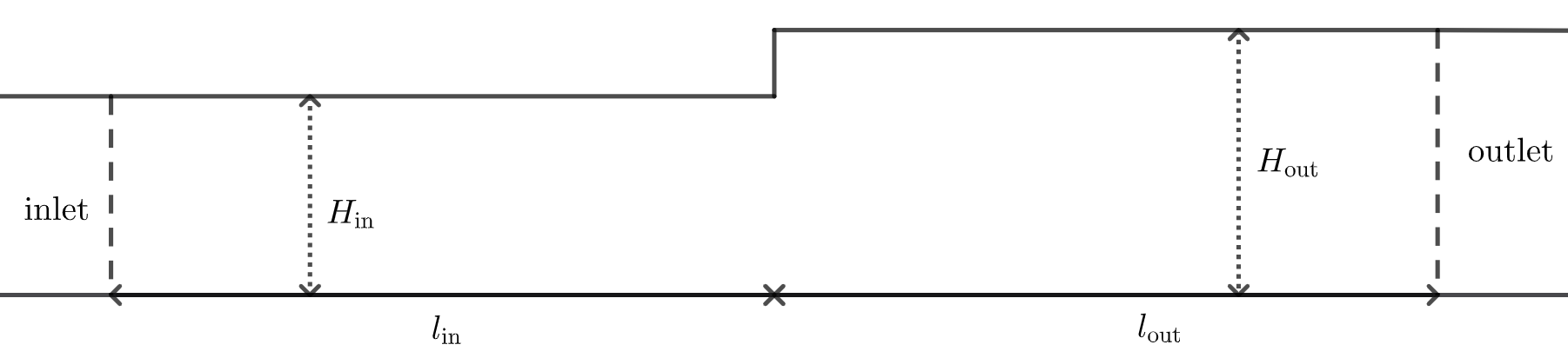}
    \caption{\color{black}A schematic of the backward facing step.}
    \label{schematic_bfs}
\end{figure}

The backward facing step is a classical example in lubrication theory; this piecewise constant height features a single discontinuity,
\begin{equation}h(x) = \begin{cases}
    H_\text{in} & 0 \le x < l_\text{in}\\
    H_\text{out} & l_\text{in} \le x \le L
\end{cases},\end{equation}
where $H_\text{out} > H_\text{in}$ and $L=l_\text{in}+l_\text{out}$. {\color{black}A schematic for the backward facing step geometry is shown in \cref{schematic_bfs}.} Both the PWC and the PWL methods for the Reynolds equation give the exact solution to the backward facing step. {\color{black}However, the FD method for the Reynolds equation suffers a loss of accuracy due to the discontinuity in the height function;  we consider convergence testing for all methods in \cref{app_convergence}.}

\begin{figure}[h]
    \centering
    \subfloat[Pressure contours]{
    \includegraphics[width=.95\textwidth]{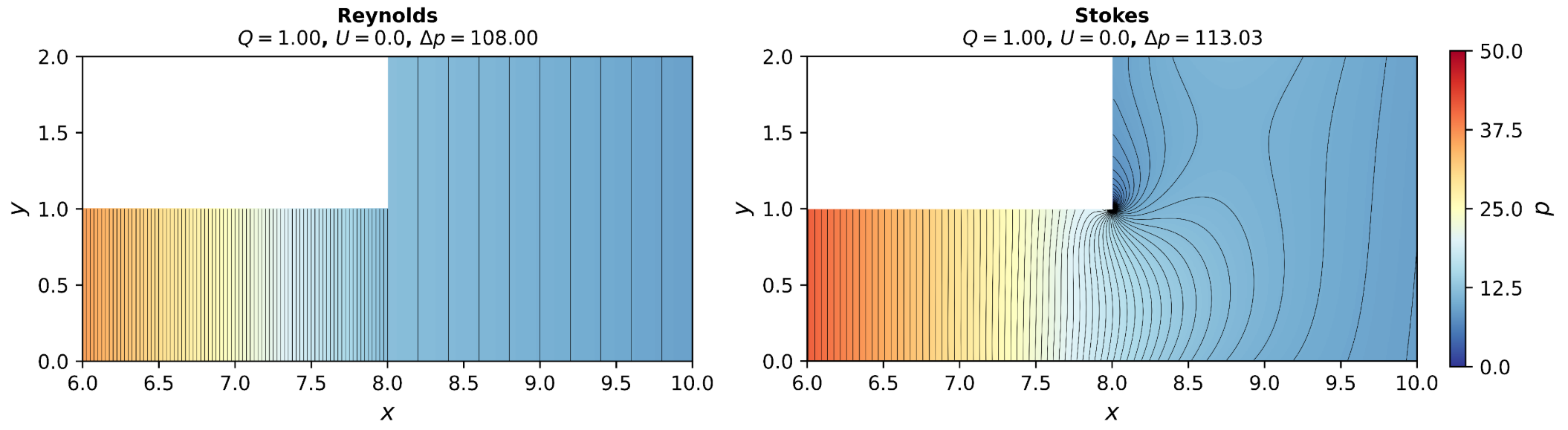}\label{bfs_p}}
    
    \subfloat[Streamlines]{
    \includegraphics[width=.95\textwidth]{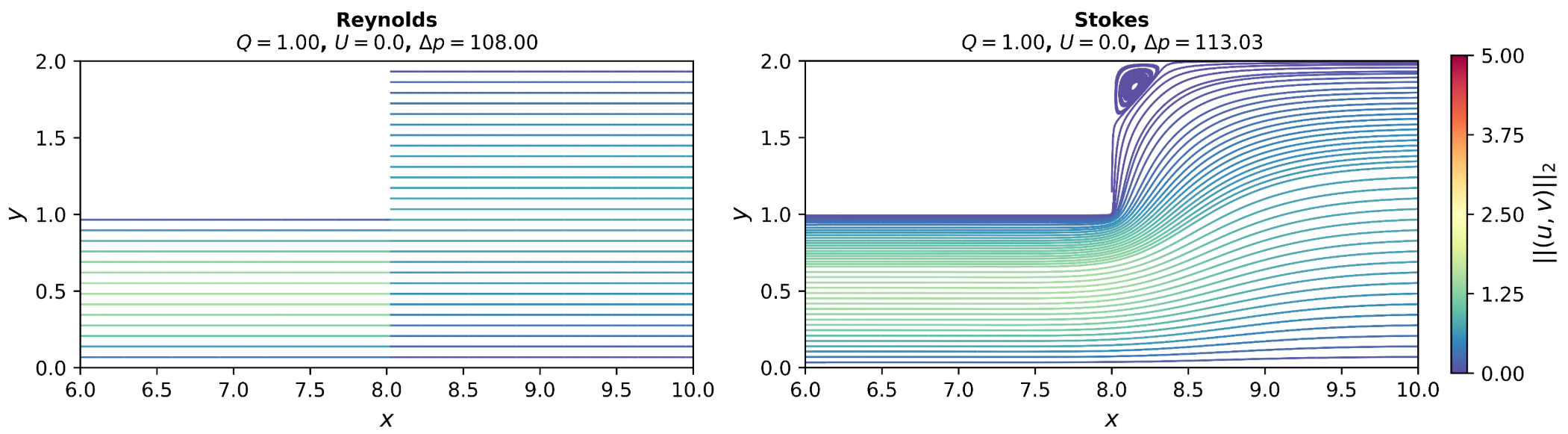}\label{bfs_v}}
    \caption{The pressure and velocity solutions for the backward facing step with the Reynolds equation (left) and the Stokes equation (right). The solution to the Reynolds equation underestimates the pressure drop $\Delta p$, and does not capture corner flow recirculation.}\label{bfs}
\end{figure}

The pressure and velocity solutions to the backward facing step with the Reynolds and the Stokes equations are shown in \cref{bfs} for $H_\text{in}=1$, $H_\text{out}=2$, $l_\text{in}=8$, $L=16$, and with the boundary conditions $\mathcal{Q} = 1$, $\mathcal{P}_N=0$, $\mathcal{U} = 0$. {\color{black}\cref{bfs_p} depicts the pressure contours, and \cref{bfs_v} depicts the streamlines colored with the magnitude of velocity $||(u,v)||_2$.} The solution to the Stokes equation displays significant cross film pressure variation in the vicinity of the step that the Reynolds equation does not capture. The solution to the Reynolds equation underestimates the total pressure drop $\Delta p$. Moreover, velocity for the Stokes equation depicts corner flow recirculation, whereas the velocity solution from the Reynolds equation is discontinuous at the step due to the surface discontinuity. 

\subsection{Piecewise linear height}
\begin{figure}
    \centering
    \includegraphics[width=0.9\textwidth]{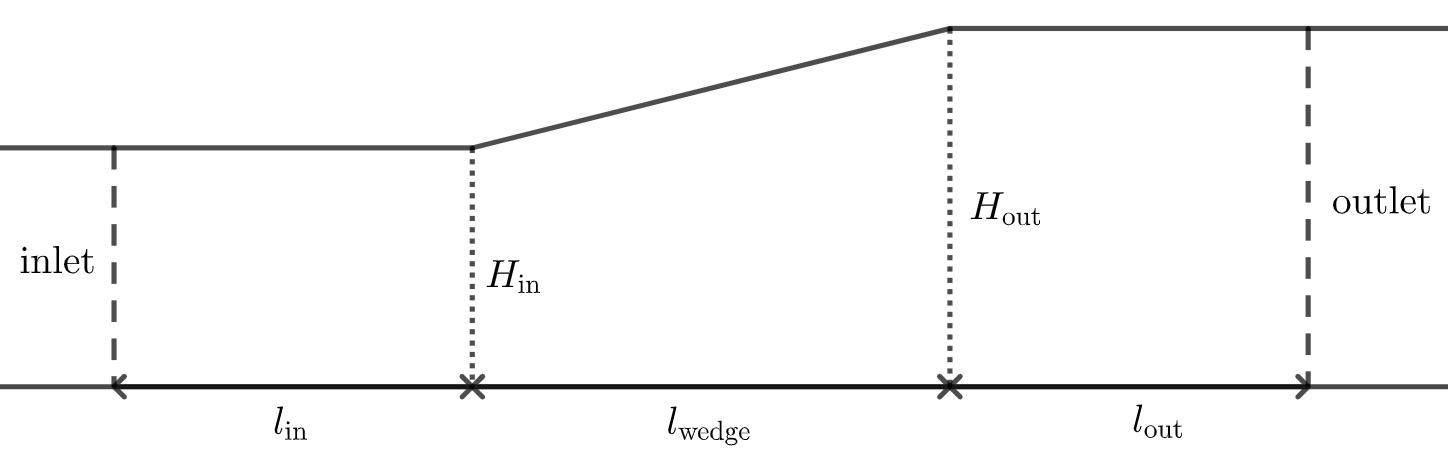}
    \caption{\color{black}A schematic of the wedge slider.}
    \label{schematic_wedge}
\end{figure}

Next we consider the wedge slider; the continuous piecewise linear height is given by,
\begin{equation}
    h(x) = \begin{cases}
    H_\text{in} & 0 \le x < l_{\text{in}}\\
     H_\text{in}+\frac{H_\text{out}-H_\text{in}}{l_\text{wedge}}(x-l_\text{in}) & l_\text{in}\le x < l_\text{in} + l_\text{wedge}\\
    H_\text{out} & l_\text{in} + l_\text{wedge} \le x \le L
\end{cases},\end{equation}
where $L=l_\text{in}+l_\text{wedge}+l_\text{out}$, $H_\text{out} > H_\text{in}$ and $l_\text{wedge}>0$. {\color{black} A schematic for the wedge slider is shown in \cref{schematic_wedge}}.
Here, the PWL method gives the exact solution to the Reynolds equation, and both the PWC method and the FD method converge to the exact solution with second-order accuracy.

\begin{figure}[h]
    \centering
    \subfloat[Pressure contours]{
    \includegraphics[width=.95\textwidth]{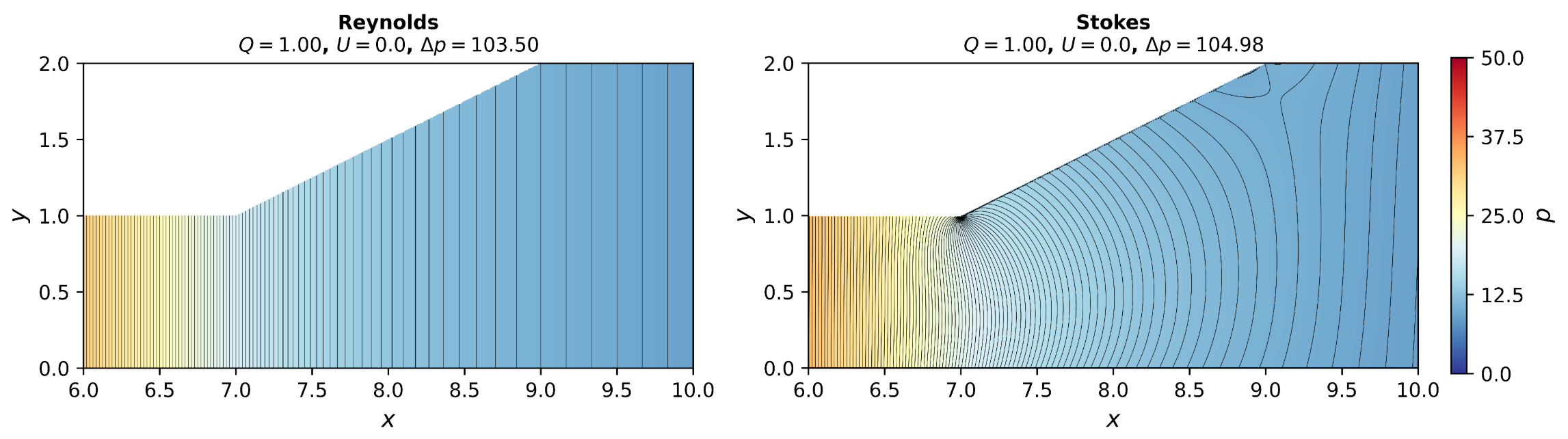}}
    
    \subfloat[Streamlines]{
    \includegraphics[width=.95\textwidth]{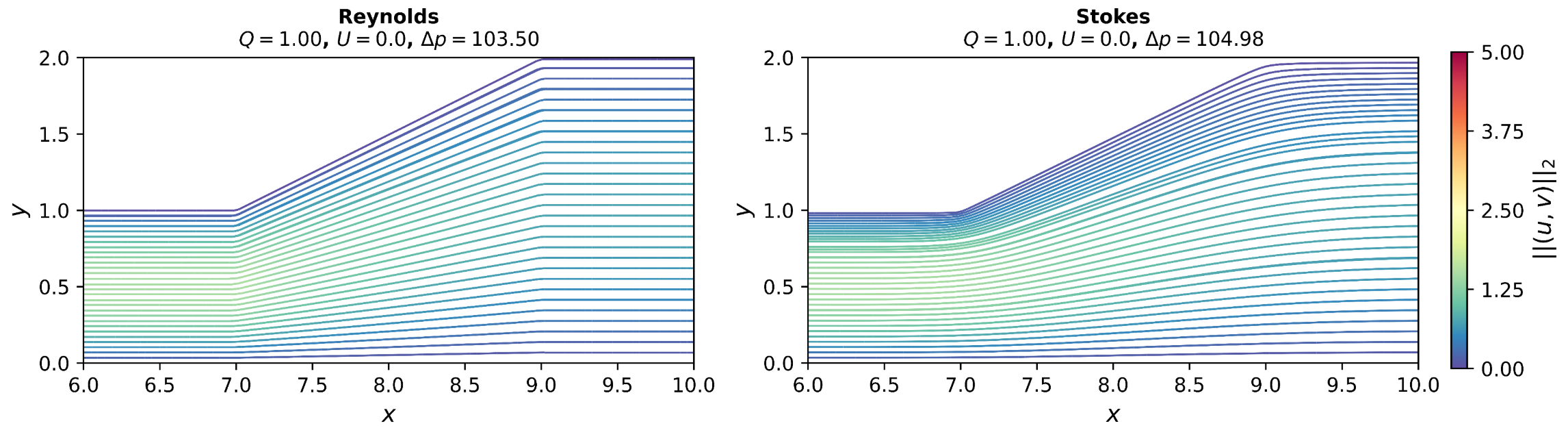}}
    \caption{The pressure and velocity solutions for the wedge slider ($l_\text{wedge}=2$)  with the Reynolds equation (left) and the Stokes equation (right). The solutions to the Reynolds and the Stokes equations are similar in the case of a wedge with moderate slope.}\label{wedge_l2}
\end{figure}

The pressure and velocity solutions to the wedge slider for the Reynolds and the Stokes equations are shown in \cref{wedge_l2} for $H_\text{in}=1$, $H_\text{out}=2$, $l_\text{wedge}=2$, $l_\text{in}=l_\text{out}$, $L=16$, and with the boundary conditions $\mathcal{Q} = 1$, $\mathcal{P}_N=0$, $\mathcal{U} = 0$. For this case of a moderately sloped wedge, the Reynolds equation is a fair approximation to the Stokes equation. However, the vertical velocity $v$ corresponding to the Reynolds equation is not $x$-continuous at the discontinuities in the surface gradient, furthermore, the discontinuity in $v$ is more pronounced when the height is larger.

\begin{figure}[h]
    \centering
    \subfloat[Pressure contours]{
    \includegraphics[width=.95\textwidth]{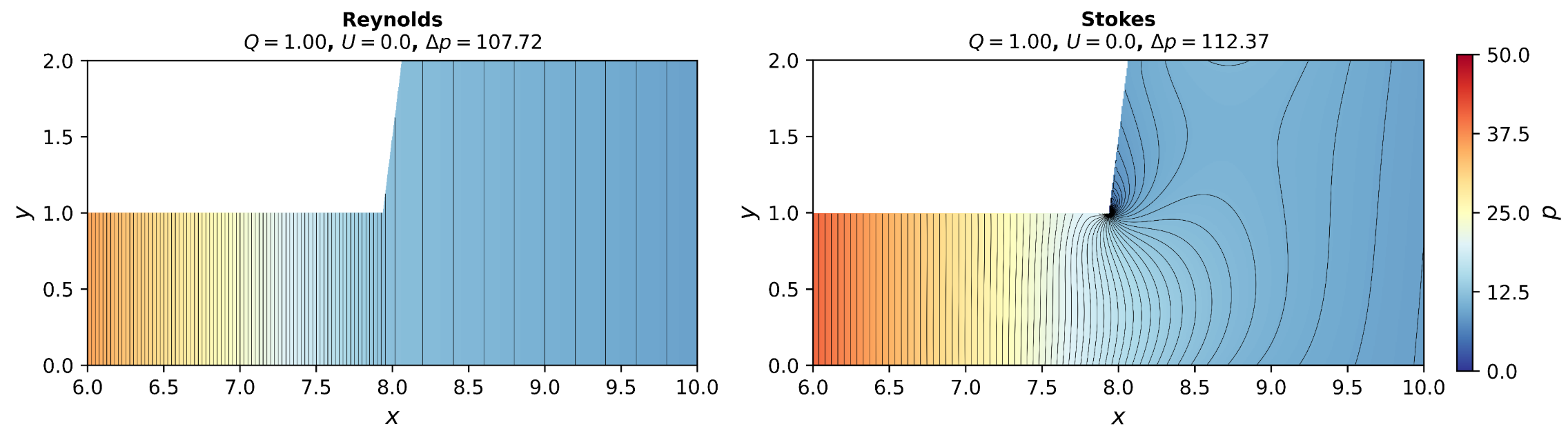}}
    
    \subfloat[Streamlines]{
    \includegraphics[width=.95\textwidth]{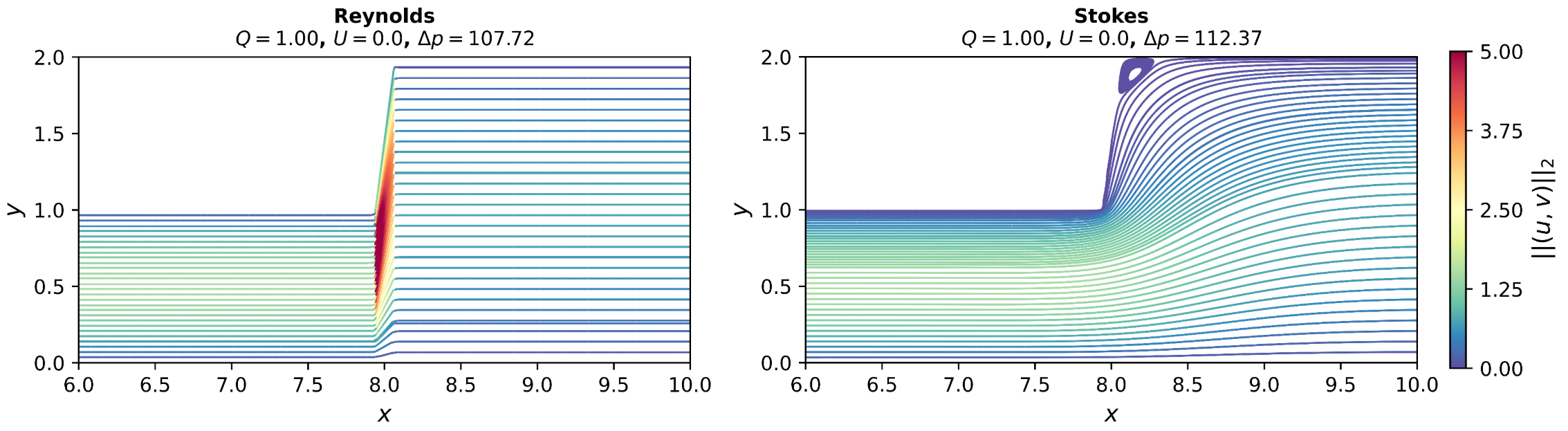}}
    \caption{\color{black}The pressure and velocity solutions for the wedge slider ($l_\text{wedge}=0.125$) with the Reynolds equation (left) and the Stokes equation (right). The solutions to the Reynolds and the Stokes equations are dissimilar in the case of a wedge with a steep slope.}\label{wedge_l1/8}
\end{figure}

{\color{black}As the slope of the wedge increases, the wedge slider approaches the solution to the backward facing step. \cref{wedge_l1/8} depicts the pressure and velocity solutions to the wedge slider for the Reynolds and the Stokes equations, now with $l_\text{wedge}=0.125$ to induce a steeper slope. The remaining parameters are consistent with the previous example; the total length is fixed at $L=16$ and $l_\text{in}=l_\text{out}=(L-l_\text{wedge})/2$. For the wedge slider with a steep slope, the solutions of the Reynolds and Stokes equations clearly differ. In the solution to the Reynolds equation, the discontinuities in $v$ are pronounced and the magnitude of velocity in the region of the large surface gradient is amplified. Moreover, the corner flow recirculation seen in the Stokes solution to the steep wedge slider and the backward facing step is not observed in the solution to the Reynolds equation.}

\subsection{Piecewise smooth heights}
Next, we consider smooth and piecewise smooth geometries, including a logistic step and a sinusoidal slider. For these examples, the PWC, PWL, and FD methods for the Reynolds equation are all numerical approximations. Convergence testing against an exact solution for the sinusoidal slider is presented in \cref{app_convergence}.

To utilize the PWL and PWC solvers for heights which are not piecewise constant or piecewise linear, we consider suitable approximations to the height function. Where the FD method uses $N+1$ grid points $\{x_i\}_{i=0}^N$ and discretizes the height as $\{h_i\}_{i=0}^N$, the PWC method considers $N$ intervals $\{[x_i,x_{i+1}]\}_{i=0}^{N-1}$ and constant heights $\{(h_{i+1}+h_i)/2\}_{i=0}^{N-1}$ for each interval. Likewise, the PWL method considers $N$ intervals $\{[x_i,x_{i+1}]\}_{i=0}^{N-1}$ with the same height discretization $\{h_i\}_{i=0}^N$ as the FD method, corresponding to $N$ constant height gradients $\{\frac{\Delta h}{\Delta x}\big|_i\}_{i=0}^{N-1}$ for each interval. 

\subsubsection{The logistic step}
\begin{figure}
    \centering
    \includegraphics[width=0.9\textwidth]{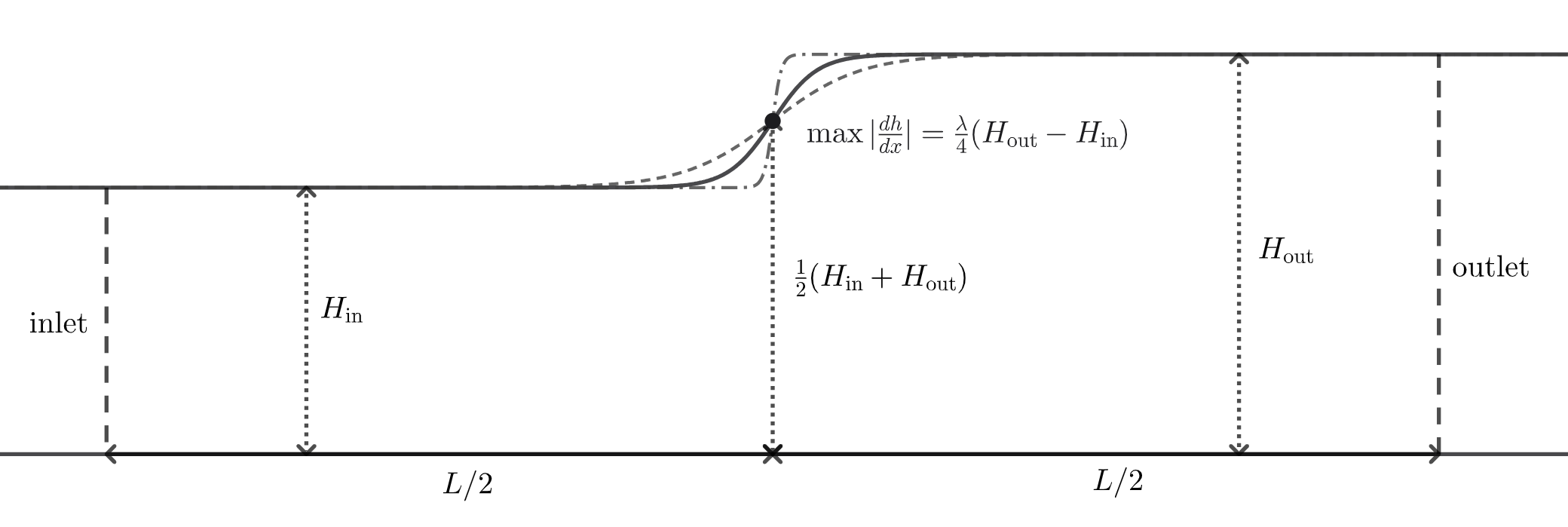}
    \caption{\color{black}A schematic of the logistic step.}
    \label{schematic_logistic}
\end{figure}

The logistic step is a smooth analogue to the backward facing step and wedge slider. The height is given by,
\begin{equation}\label{logistic_eq} h(x) = H_{\text{in}} + \frac{H_{\text{out}}-H_{\text{in}}}{1+e^{\lambda \big(L/2-x\big)}},\end{equation} where $H_\text{out} > H_\text{in}>0$ are the inlet and outlet heights, and $\lambda$ corresponds to the maximum surface gradient, $\max |\frac{dh}{dx}| = \lambda(H_\text{out}-H_\text{in})/4$ at the midpoint $x=L/2$. {\color{black}A schematic for the logistic step is shown in \cref{schematic_logistic}.} {\color{black} As $\lambda\to0$, the surface gradient decreases in magnitude but the range in $x$ of significant surface variation increases. As $\lambda\to\infty$, the logistic step approaches the discontinuous backward facing step. Note that the length $L$ is chosen large enough relative to $\lambda$ and $|H_\text{out}-H_\text{in}|$ such that the inlet and outlet heights are achieved to a tolerance much smaller than the discretization error, i.e. $|h(x_0)-H_\text{in}|=|h(x_L)-H_\text{out}|\ll \Delta x/2$.} 

\begin{figure}[h]
    \centering
    \subfloat[Pressure contours]{
    \includegraphics[width=.95\textwidth]{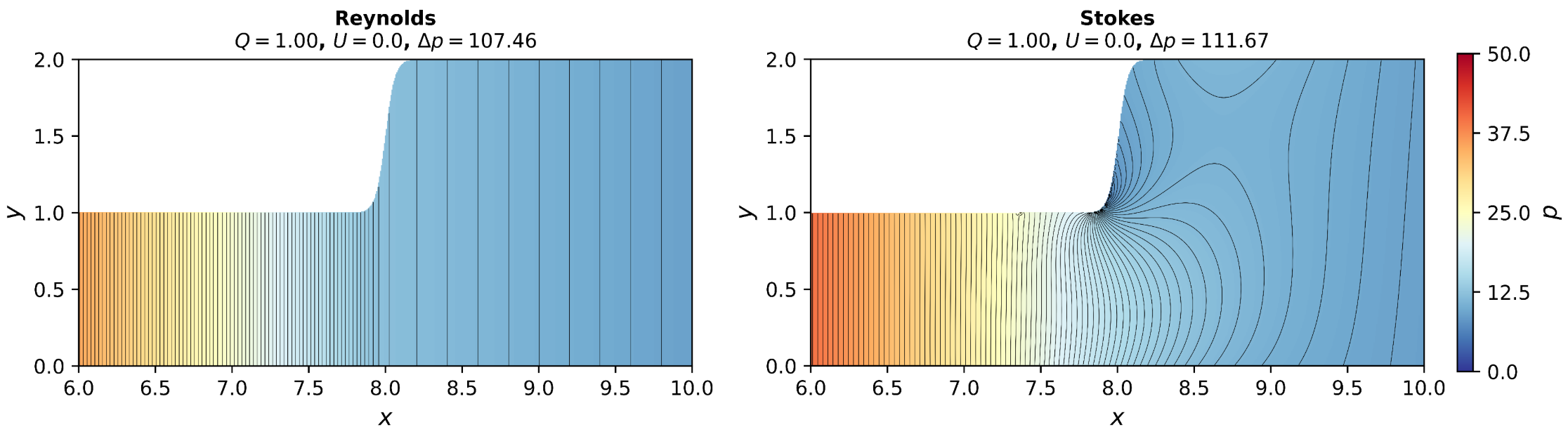}}
    
    \subfloat[Streamlines]{
    \includegraphics[width=.95\textwidth]{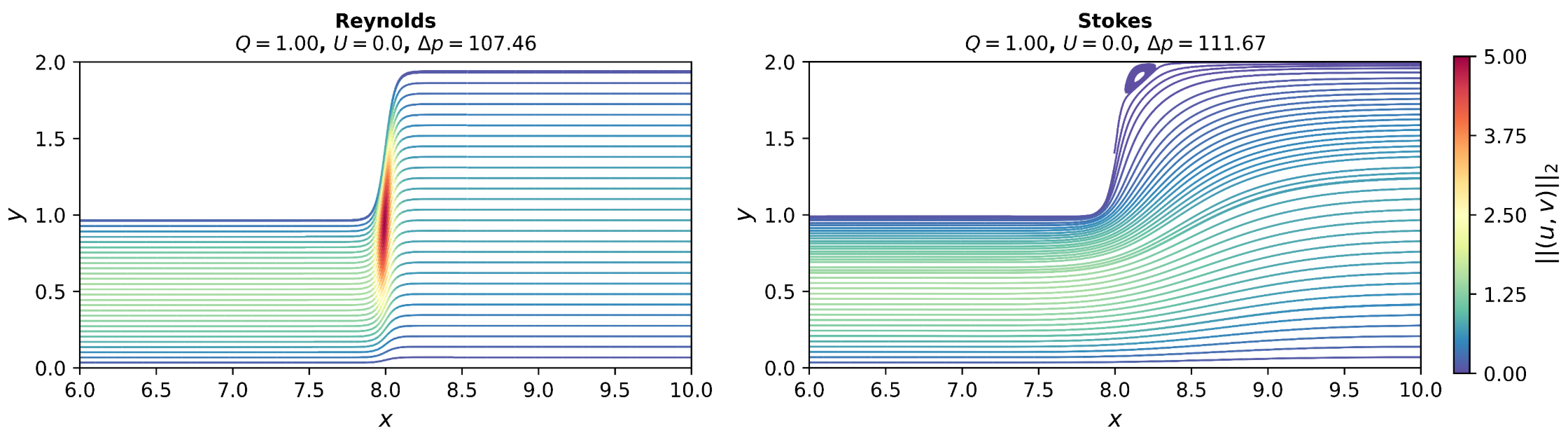}}
    \caption{The pressure and velocity solutions for the logistic step with the Reynolds equation (left) and the Stokes equation (right). As with the backward facing step, the Reynolds equation underestimates the pressure drop $\Delta p$, and does not capture cross film pressure variation or corner flow recirculation as seen with the Stokes equation.}\label{logistic}
\end{figure}

The pressure and velocity solutions to the logistic step for the Reynolds and the Stokes equations are shown in \cref{logistic} for $H_\text{in}=1$, $H_\text{out}=2$, $\lambda=32$, $L=16$, and with the boundary conditions $\mathcal{Q} = 1$, $\mathcal{P}_N=0$, $\mathcal{U} = 0$. {\color{black}As with the wedge slider, when the magnitude of surface variation is increased, or when the range in $x$ of the surface variation is decreased, the thin film assumptions of lubrication theory break down and the Reynolds and Stokes solutions begin to disagree.} The pressure for the Stokes equation has significant cross film pressure variation $\frac{\partial p}{\partial y}$ in the vicinity of the large surface gradient, whereas the pressure for the Reynolds equation is necessarily one dimensional. At this steep of slope, the velocity for the Stokes equation depicts flow recirculation similar to the backward facing step in \cref{bfs}. The velocity for the Reynolds equation does not have flow recirculation, and the velocity magnitude is overestimated in the vicinity of the large surface gradient. 

\subsubsection{The sinusoidal slider}
\begin{figure}
    \centering
    \includegraphics[width=0.9\textwidth]{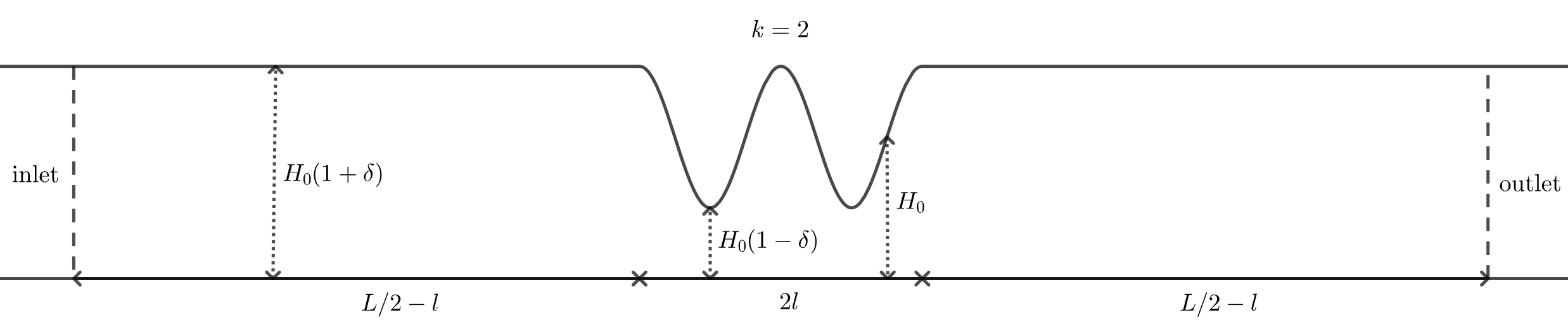}\\\vspace{.5em}
    \includegraphics[width=0.9\textwidth]{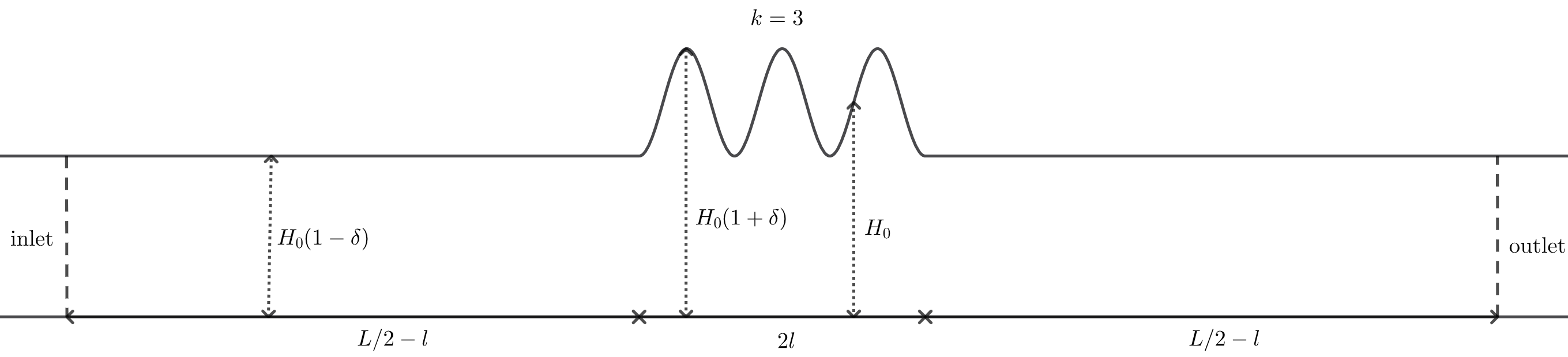}
    \caption{\color{black}Schematics of the sinusoidal slider with $k=2$ (upper) and $k=3$ (lower).}
    \label{schematic_sinusoid}
\end{figure}

Compared with the previous examples which are all variations on a step texture, the sinusoidal slider exhibits a cavity texture. The height is,
\begin{equation}
    h(x) = \begin{cases}
        H_0 (1 + \delta) & l<|x| <L/2,  \text{ $k$ even}\\
        H_0 (1 - \delta) & l<|x| <L/2,  \text{ $k$ odd}\\
        H_0 (1 + \delta \cos(x \pi k/l)) & |x|\le l
    \end{cases},
\end{equation}
where $k$ is the integer wave number on length $2l$, $H_0$ is the equilibrium height, $\delta\in[0,1)$ gives the amplitude $\delta H_0$, and $L$ is the total length. {\color{black}A schematic of the sinusoidal slider is shown in \cref{schematic_sinusoid}.} The regions of constant height at the inlet and outlet are aligned with the extrema of the sinusoid, keeping the surface gradient continuous. When $k$ is odd, $k$ is the number of positive textures and $\min h(x)=H_0 (1 - \delta)$ is the inlet and outlet height. When $k$ is even, $k$ is the number of negative textures and $\max h(x)=H_0(1+\delta)$ is the inlet and outlet height.  

\begin{figure}[h]
    \centering
    \subfloat[Pressure contours]{
    \includegraphics[width=.95\textwidth]{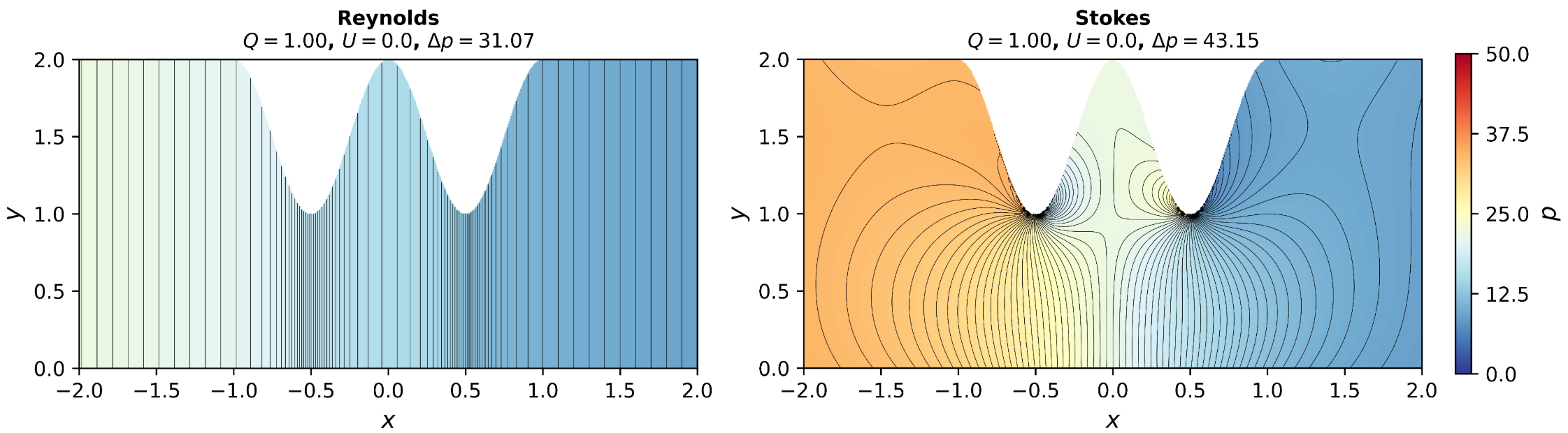}}
    
    \subfloat[Streamlines]{
    \includegraphics[width=.95\textwidth]{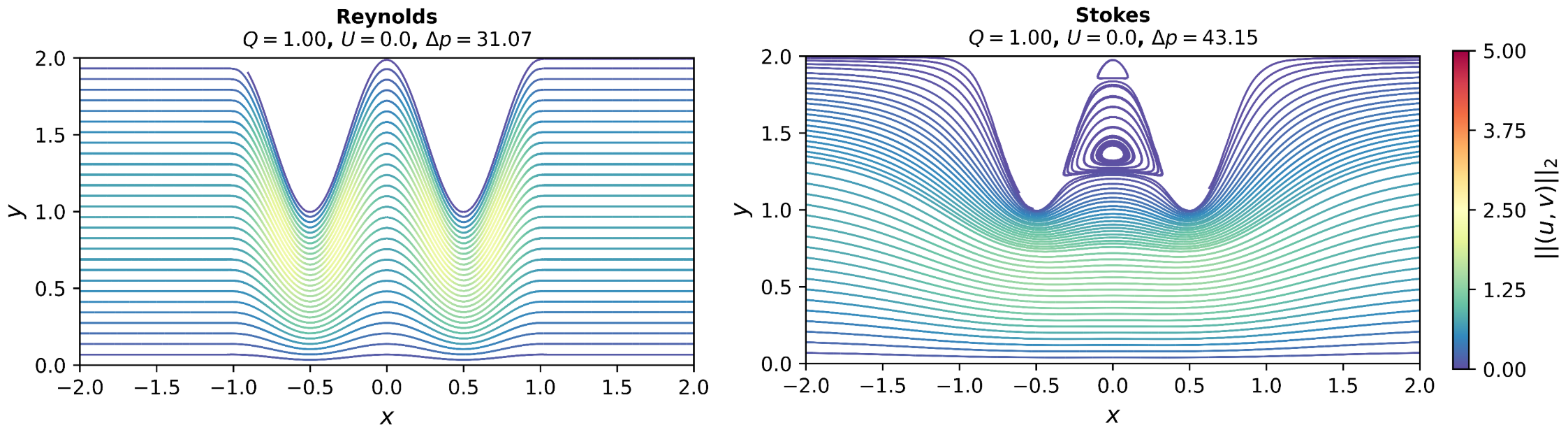}}
    \caption{The pressure and velocity solutions for the sinusoidal slider with the Reynolds equation (left) and the Stokes equation (right). The Reynolds solution overestimates the velocity magnitude when the surface gradient is large relative to the height.}\label{sinusoid}
\end{figure}

The pressure and velocity solutions to the sinusoidal slider for the Reynolds and the Stokes equations are shown in \cref{sinusoid} for $H_0=3/2$, $\delta = 1/3$, $k =2$, $l=1$, $L=16$, and with the boundary conditions $\mathcal{Q} = 1$, $\mathcal{P}_N=0$, $\mathcal{U} = 0$. The velocity for the Stokes equation depicts primary and secondary flow recirculation in the region of positive texturing that the Reynolds equation does not capture. The velocity magnitude for to the Stokes equation is largest at the valleys of the sinusoid, whereas for the Reynolds equation, the velocity magnitude is largest either side of these valleys where the surface gradient is large. The overall pressure drop $\Delta p$ of the Reynolds solution is significantly smaller than that of the Stokes solution. In general, for regions where the surface gradient is large relative to the height, the solution to the Stokes equation has significant cross film pressure variation $\frac{\partial p}{\partial y}$ and the Reynolds equation tends to overestimate the velocity magnitude.

\section{Timing}
When the height is piecewise constant or piecewise linear, the PWC and PWL methods give the exact solution to the Reynolds equation and are certainly faster than the Reynolds FD solution. The PWC and PWL methods scale with the number of piecewise components, whereas the FD method requires a fine grid resolution to capture steep slopes and surface discontinuities.

\begin{figure}[h]
    \centering
    \includegraphics[width=.95\textwidth]{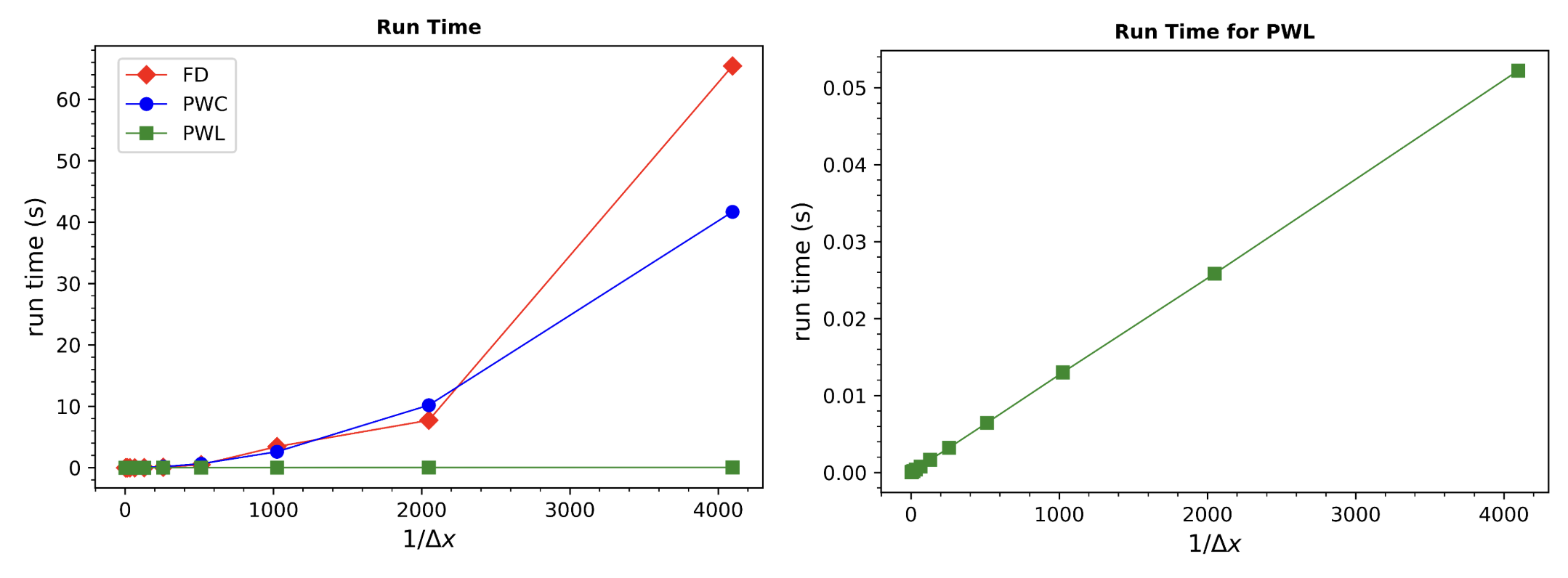}
    \caption{The PWL method is linear time and faster than the PWC method and the FD method.}\label{runtimes}
\end{figure}

When the height is non-linear, the PWC, PWL and FD methods all approximate the solution on a grid. Recall, where the FD method uses $N+1$ grid points, the PWC and PWL methods use $N$ piecewise components; the grid spacing in both cases is $\Delta x = 1/N$. \cref{runtimes} shows the computational run times for the FD, PWC, and PWL methods with the logistic step example. For a non-linear height, the PWL method is significantly faster than both the PWC and FD methods. The FD method is evaluated with numpy.linalg.solve, which is an LU method with partial pivoting, having $\mathcal{O}(N^3)$ time complexity. In the PWC method, a nested loop is used to compute each element of the Schur complement inverse $K^{-1}_{i,j}$ when we evaluate $M^{-1}{\bf b}$, leading to $\mathcal{O}(N^2)$ complexity. In the PWL method, the upper triangular Schur complement is simple invert and $M^{-1}{\bf b}$ can be evaluated in a single loop, corresponding to $\mathcal{O}(N)$ complexity. These time complexities are confirmed by the results in \cref{runtimes}. {\color{black} Note that we have experimented with various other iterative solvers for the Reynolds finite difference formulation and found the performance to be comparable to the LU solver. The PWL approach has the advantage of having an explicit inverse, so that compared with the FD method, the PWL method requires no matrix solve.}

\section{Conclusions}
We presented the piecewise constant (PWC) and piecewise linear (PWL) methods for the Reynolds equation and compared with the finite difference (FD) method. The PWC and PWL methods take advantage of the exact solution to the Reynolds equation when the height is piecewise linear, and couple the exact solution for each piecewise component using the assumption of constant flux and continuous pressure. The PWC and PWL methods seamlessly handle discontinuities in the surface and retain their second-order accuracy, compared with the FD method which has reduced order of accuracy in the case of surface discontinuities. The PWL method takes advantage of the constant flux to model the pressure gradients on each piecewise component, whereas the PWC method treats each constant pressure gradient as a solution variable. As a result, the PWL method is the most computationally efficient, running with linear time complexity. 

We then compared the solutions of pressure and velocity for the Reynolds and the Stokes equations. The pressure from the Stokes equation exhibits significant cross film pressure variation in the presence of large surface gradients or surface discontinuities. The pressure from the Reynolds equation has no cross film pressure variation, leading to an underestimation of the overall pressure drop in the case of large surface gradients or surface discontinuities. For the backward facing step, and for the logistic step, wedge slider and sinusoidal slider with steep slopes, the velocity from the Stokes equation depicts corner flow recirculation. The velocity from the Reynolds equation does not correctly model flow recirculation compared with the Stokes equations, and the magnitude of velocity from the Reynolds equation is overestimated when large surface gradients occur alongside small heights. Moreover, the velocities from the Reynolds equation are discontinuous at discontinuities in the surface, and the vertical velocity is also discontinuous at discontinuities in the surface gradient. Ultimately, the Reynolds equation is best suited for the case of continuous and slowly varying surface geometries. 

\section*{Acknowledgments}
This work was supported by National Science Foundation (NSF) grant DMS-2512565 to TGF. We also acknowledge use of the Brandeis High Performance Computing Cluster (HPCC) which is partially supported by the NSF through DMR-MRSEC 2011846 and OAC-1920147.

\appendix

\section{Finite difference method for the Reynolds equation}\label{app_reyn_fd} Here we present a finite difference method (FD) for the Reynolds equation.
Define the uniform discretisation of the domain $[x_0, x_L]$,
\begin{align}\{x_i\}_{i=0}^N, && x_i = x_0 + i \Delta x, && \Delta x = |x_L-x_0|/N.
\end{align}
A second-order accurate difference approximation for the Reynolds equation \cref{reynolds} is,
 \begin{multline}\label{reyn_fd}
\Big(h_{i+1}^3 + h_i^3\Big)p_{i+1} - \Big(h_{i+1}^3 + 2h_i^3 +h_{i-1}^3\Big)p_i + \Big(h_i^3 +h_{i-1}^3\Big)p_{i-1} \\= 6 \eta \mathcal{U}\Big(h_{i+1}-h_{i-1}\Big)\Delta x,
 \end{multline}
 where $h_i = h(x_i)$ and $p_i = p(x_i)$.
 The form \cref{reyn_fd} is obtained by approximating the outer derivative on the left-hand-side of \cref{reynolds} on half grid points,
 \begin{equation}\Bigg[\frac{d}{dx}\bigg[\big[h(x)\big]^3\frac{dp}{dx}\bigg]\Bigg]_{i} = \frac{\bigg[\big[h(x)\big]^3\frac{dp}{dx}\bigg]_{i+1/2} - \bigg[\big[h(x)\big]^3\frac{dp}{dx}\bigg]_{i-1/2}}{\Delta x},\end{equation}
and using an average to approximate the height function at the half grid points,
 \begin{equation}\bigg[\big[h(x)\big]^3\frac{dp}{dx}\bigg]_{i\pm1/2} = \frac{\pm\Big(h_{i\pm1}^3+h_{i}^3\Big)\Big(p_{i\pm1}-p_{i}\Big)}{2\Delta x}.\end{equation} 
The prescribed flux $\mathcal{Q}$ and boundary velocity $\mathcal{U}$ dictate the inlet pressure gradient, incorporated through a right-sided difference, \begin{equation}\frac{dp}{dx}\bigg|_{0}=\frac{3p_0 - 4p_{1}+p_{2}}{2\Delta x}= \frac{-12\eta \mathcal{Q}}{h_0^3} + \frac{6\eta \mathcal{U}}{h_0^2}.\end{equation} The prescribed outlet pressure is simply applied as $p_N = \mathcal{P}_N$.
 
The equation \cref{reyn_fd} for $0 < i < N$ and the two boundary conditions associated with $i=0$ and $i=N$ characterize a size $N+1$ linear system solving the Reynolds equation for $\{p_k\}_0^N$. We utilize numpy.linalg.solve to evaluate this linear system.
 
\section{Finite difference method for the Stokes equation} \label{app_stokes}
The biharmonic formulation of the incompressible Navier-Stokes equations provides an effective method of solution for low Reynolds number steady state flows \cite{gupta_new_2005,marner_potential-velocity_2014,sen_4oec_2015,biswas_hoc_2017}, and results in exactly the Stokes equations at zero Reynolds number. 
Through introduction of the {\color{black}streamfunction} $\psi(x,y)$ satisfying,
\begin{align} \label{stream}
 \hfill&& u = \frac{\partial\psi}{\partial y}, \hspace{3em} v = -\frac{\partial \psi}{\partial x},&&\hfill
\end{align}
the Navier-Stokes equations,
\begin{align}\label{n-s_x}
 \frac{\partial p}{\partial x} &= \eta \Big(\frac{\partial^2 u}{\partial x^2} + \frac{\partial^2 u}{\partial y^2}\Big) - \rho\Big( \frac{\partial u}{\partial t}+u\frac{\partial u}{\partial x} + v \frac{\partial u}{\partial y}\Big),\\
 \label{n-s_y} \frac{\partial p}{\partial y} &= \eta \Big(\frac{\partial^2 v}{\partial x^2} + \frac{\partial^2 v}{\partial y^2}\Big) - \rho\Big( \frac{\partial v}{\partial t}+u\frac{\partial v}{\partial x} + v \frac{\partial v}{\partial y}\Big),
\end{align}
are expressed in dimensionless terms as,
\begin{equation}\label{n-s_psivel}
 \nabla^4\Psi = \text{Re}\big(V\nabla^2 U - U \nabla^2 V),\end{equation} where $X = x/L_x$, $Y=y/L_x$, $U=u/U_*$, $V = v/U_*$, $P = p/P_*$, $\Psi = \psi/\mathcal{Q}$, and the flow is assumed to be steady state. Note, the characteristic pressure is now $P_* = \eta U_*/L_x$ which differs from that of lubrication theory by a factor of $\varepsilon^2$.  
When $\text{Re} = 0$, the velocity-{\color{black}streamfunction} formulation \cref{n-s_psivel} reduces to the Stokes equation $\nabla^4 \psi=0$, for which it is equivalent to use dimensional variables. 

To obtain the numerical solution to the Stokes equations, we assume the inlet and outlet flow profiles correspond to a fully developed laminar flow with flux $\mathcal{Q}$, and restrict to examples with a zero height gradient in the vicinity of the inlet and outlet. The surface boundary conditions for velocity are \cref{reyn_bc_u,reyn_bc_v}, and the pressure satisfies \cref{reyn_bc_p_mixed}. The inlet and outlet velocity profiles are,
\begin{align}\label{n-s_bc_vel}
 u(x_0,y) &= u_\text{Re}(x_0,y), &&\frac{\partial u}{\partial x}\Big|_{x_L,y} = 0,\\
 v(x_0,y) &= 0, && v(x_L,y) = 0,
\end{align}
where $u_\text{Re}(x,y)$ is the Reynolds equation velocity \cref{reyn_u}  expressed in terms of $\mathcal{Q}$ from \cref{flux}. The corresponding boundary conditions for the {\color{black}streamfunction} $\psi$ are,
\begin{align}\label{n-s_bc_stream}
 \psi(x_0,y) &= \int_{0}^y u_\text{Re}(x_0,\hat{y})d\hat{y}, &&\frac{\partial \psi}{\partial x}\Big|_{x_L,y} = 0,\\
 \psi(x, 0) &= 0, & &\psi(x,h(x)) = \mathcal{Q},
\end{align}
where,
\begin{equation} \label{n-s_bc_stream_x0}
\int_{0}^y u_\text{Re}(x,\hat{y})d\hat{y} =\frac{\mathcal{Q}y^2}{\big[h(x)\big]^3}\Big(3h(x)-2y\Big) +  \frac{\mathcal{U}y}{\big[h(x)\big]^2}\Big(h(x)-y\Big)^2.
\end{equation}

The solution $\psi$, $u$, $v$ to the Stokes equation is determined through an iterative second-order accurate finite difference method. The method presented here for $\text{Re}=0$ is adapted from the method given in \cite{biswas_hoc_2017, gupta_new_2005} for the full biharmonic Navier-Stokes equations.
The Stokes equation is discretized as,
\begin{multline}\label{n-s_disc}
    28\psi_{i,j} -  8  \Big(\psi_{i-1,j}+\psi_{i+1,j}+\psi_{i,j-1}+\psi_{i,j+1}\Big) \\ + \Big(\psi_{i-1,j-1}+\psi_{i-1,j+1}+\psi_{i+1,j-1}+\psi_{i+1,j+1}\Big)
     \\ = 3\Delta x \big(u_{i,j-1} - u_{i,j+1} + v_{i+1,j}-v_{i-1,j}\big).
\end{multline}
The stream-velocity relation \eqref{stream} is discretised as,
\begin{align}\label{n-s-vel_disc_u}
    u_{i,j} &= \frac{-3}{4\Delta x}\big(\psi_{i,j-1}-\psi_{i,j+1}\big) - \frac{1}{4}\big(u_{i,j-1}+u_{i,j+1}\big) 
    \\ \label{n-s-vel_disc_v}
    v_{i,j} &= \frac{3}{4\Delta x}\big(\psi_{i-1,j}-\psi_{i+1,j}\big)-\frac{1}{4}\big(v_{i-1,j}+v_{i+1,j}\big).
\end{align}
In \cite{dennis_separation_2026}, we present a method for approximation of the {\color{black} streamfunction} and velocity boundary conditions with non-rectilinear domains. 

Once the {\color{black} streamfunction} and velocity solution has sufficiently converged, the pressure partial derivatives are determined using a centered second-order accurate finite difference discretisation of the Navier-Stokes equations \cref{n-s_x,n-s_y}. For \cref{n-s_x},
\begin{multline}\label{n-s-x_disc}\frac{\partial p}{\partial x}\Big|_{i,j} = \frac{\eta}{\Delta x^2} \Big(u_{i-1,j} + u_{i+1,j} - 4u_{i,j} + u_{i,j-1} + u_{i,j+1}\Big), \\ -\frac{\rho}{2\Delta x}\Big(u_{i,j}(u_{i+1,j}-u_{i-1,j}) + v_{i,j}(u_{i,j+1}-u_{i,j-1})\Big),\end{multline} and $\frac{\partial p}{\partial y}$ as in \cref{n-s_y} is similar. The pressure $p(x,y)$ is then determined by numerical integration, using the outlet boundary condition $p(x_L, 0) = \mathcal{P}_N$. Note, the pressure drop for the Stokes solution is evaluated at the lower surface, $\Delta p = p(x_0,0) - p(x_L,0)$. Since we assume a fully developed laminar flow at the inlet and outlet, the cross film pressure variation at the inlet and outlet is negligible.

\section{Convergence}\label{app_convergence}
{\color{black} Here we consider the convergence of the FD, PWL and PWC methods for the Reynolds equation, and of the iterative finite difference method for the Stokes equation.}

\subsection{Methods for the Reynolds equation}
\begin{figure}[h]
    \centering
    \includegraphics[width=.95\textwidth]{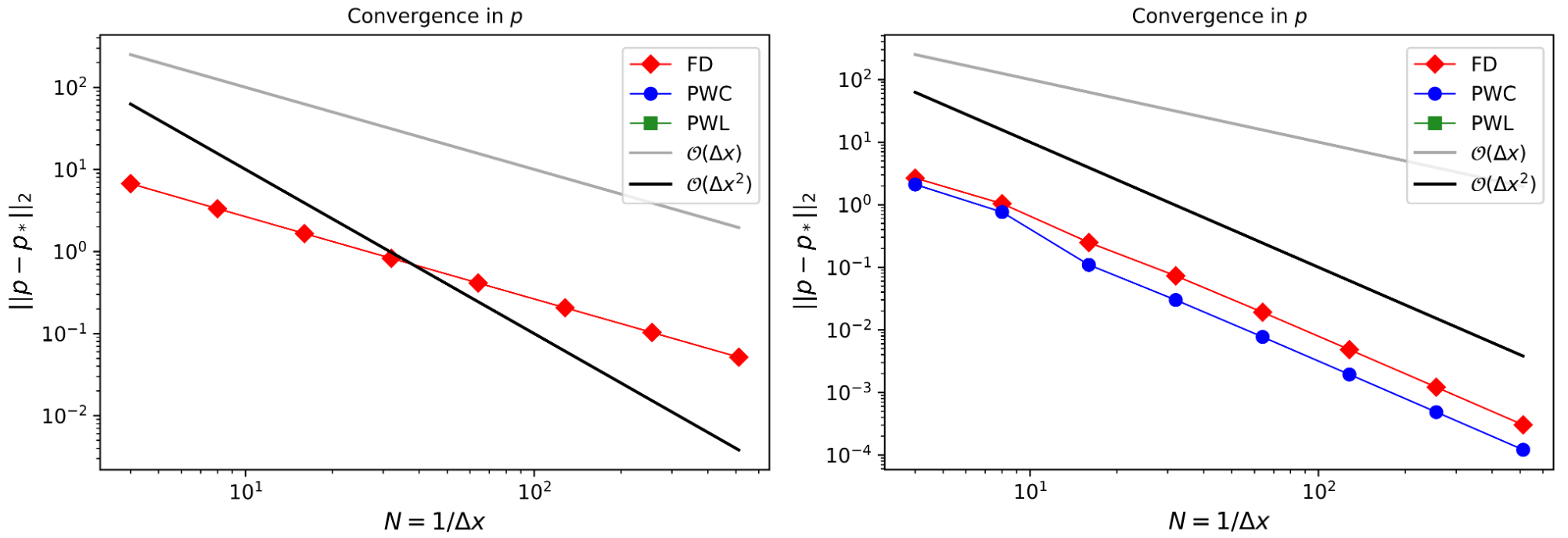}
    \caption{\color{black}Convergence of solutions to the Reynolds equation for the backward facing step (left) and the wedge slider (right). The FD method is first order accurate for the backward facing step, and second order accurate for the wedge slider. The PWC method is exact for the backward facing step and is second order accurate for the wedge slider.}\label{convergence_piecewise}
\end{figure}

{\color{black}Recall that for the backward facing step, the PWL and PWC methods give the exact solution to the Reynolds equation; likewise, for the wedge slider, the PWL method gives the exact solution to the Reynolds equation. The convergence results of the FD method applied to the backward facing step and of the FD and PWC methods applied to the wedge slider are shown in \cref{convergence_piecewise}. For the backward facing step, the FD method demonstrates first-order convergence in $p$ in the $L^1$, $L^2$ and $L^\infty$ norms. For the wedge slider, both the FD and PWC methods demonstrate second-order convergence in $p$ in the $L^1$, $L^2$ and $L^\infty$ norms. We note that although our finite difference method only obtains its full accuracy in smooth geometries, previous work has shown that local mesh refinement allows for higher-order accuracy even in the presence of discontinuities \cite{he_modified_2023}.}

To confirm that the PWL method for the Reynolds equation is also second-order accurate for arbitrary domains, we consider the sinusoidal slider for which the Reynolds equation has an exact solution \cite{takeuchi_extended_2019}. The height is given by,
\begin{equation}
    h(x) = H_0 (1 + \delta \cos(\alpha x))
\end{equation}
where $H_0>0$ is the equilibrium height, $H_0(1+\delta)$ is the maximum height ($\delta > 0$), and $\alpha\ne 0$ is the period.  
The exact pressure for the Reynolds equation assumes a fixed pressure boundary condition $\Delta p = 0$. The pressure is given by,
\begin{equation}
    p(x) = \frac{-6\eta\mathcal{U}\delta (H_0 + h(x))\sin(\alpha x)}{\alpha H_0^3(2+\delta^2)(1+\delta \cos(\alpha x))^2}
\end{equation}
From \cref{flux} or \cref{dP}, the flux $\mathcal{Q}$ corresponding to $\Delta p = 0$ is $\mathcal{Q} = \frac{\mathcal{U}\eta H_0(1-\delta^2)}{2+\delta^2}$.

\begin{figure}[h]
    \centering
    \includegraphics[width=.95\textwidth]{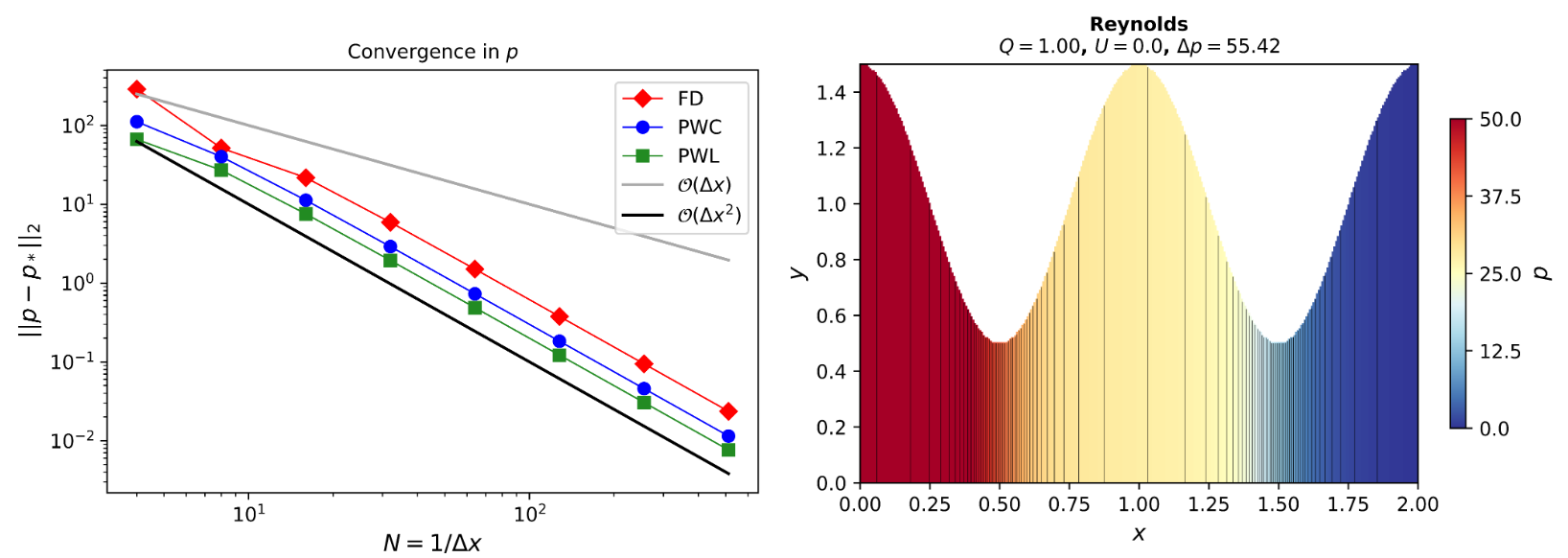}
    \caption{Convergence of solutions to the Reynolds equation for the sinusoidal slider. The FD, PWC and PWL methods all demonstrate second-order accuracy.}\label{convergence_sinusoid}
\end{figure}

The convergence tests for the sinusoidal slider confirm that the three methods for the Reynolds equation have second-order convergence in $p$ with $L^1$, $L^2$, and $L^\infty$ norms. \cref{convergence_sinusoid} shows the rate of convergence in $p$ with $L^2$ norm for the sinusoidal slider with $H_0 = 1$, $\alpha=2\pi$, $\delta=1/2$, $\mathcal{U}=0$, $\mathcal{Q} = 1$.

\subsection{Method for the Stokes equation}
{\color{black}In the absence of an exact solution for the Stokes equation, we consider convergence of the iterative finite difference method measured against the highest available grid resolution. The convergence results for the backward facing step and for the logistic step with the finest grid $N=1/\Delta x =160$ are shown in Figure C3. The streamfunction is shown to converge with a convergence rate between first and second-order over the grid resolutions tested; higher grid resolutions would be desirable for further testing, however computational resources are a limitation.}

\begin{figure}[h]
    \centering
    \includegraphics[width=.95\textwidth]{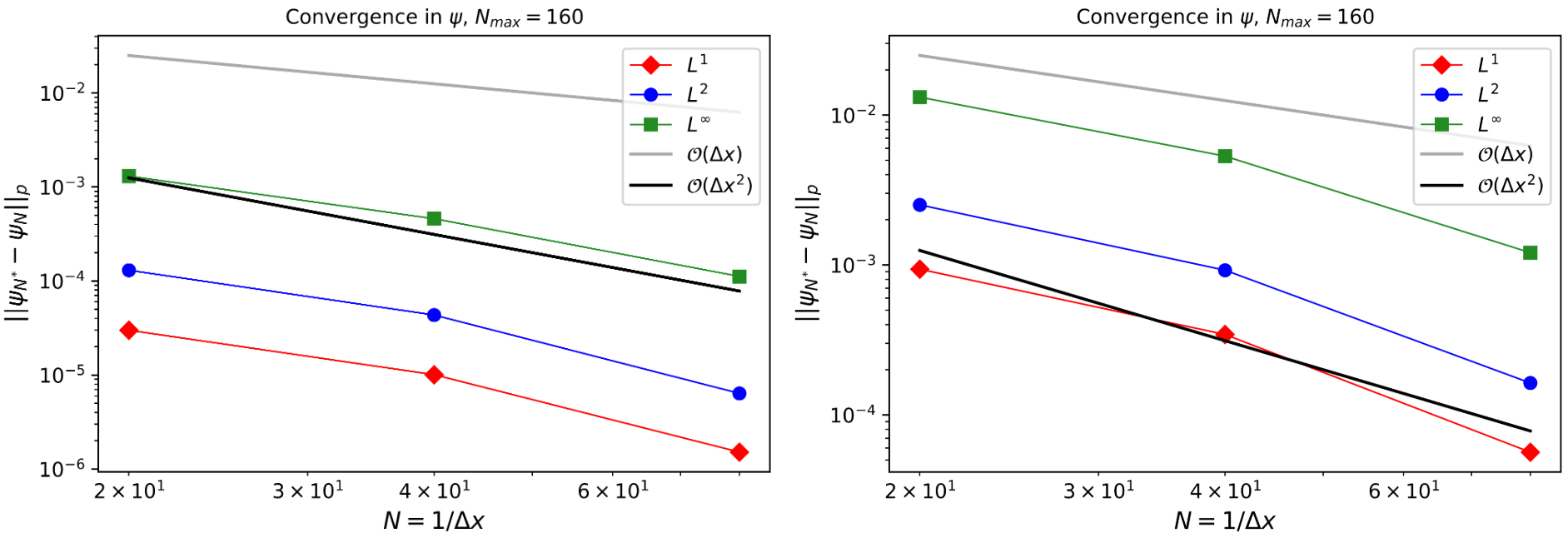}
    \caption{{\color{black}Convergence of the Stokes equation finite difference method for the backward facing step (left) and the sinusoidal slider (right). The streamfunction demonstrates a convergence rate between first and second-order over the grid resolutions tested.}}\label{convergence_stokes}
\end{figure}
\bibliography{ReynoldsStokes_bib}

@article{shyu_numerical_2018,
	title = {A numerical study on the negligibility of cross-film pressure variation in infinitely wide plane slider bearing, {Rayleigh} step bearing and micro-grooved parallel slider bearing},
	volume = {137},
	doi = {10.1016/j.ijmecsci.2018.01.031},
	journal = {Mechanical Sciences},
	author = {Shyu, Shiuh-Hwa and Hsu, Wei-Chun},
	month = mar,
	year = {2018},
	pages = {315--323}
}

@article{dobrica_about_2009,
	title = {About the validity of {Reynolds} equation and inertia effects in textured sliders of infinite width},
	volume = {223},
	doi = {10.1243/13506501JET433},
	journal = {Proceedings of the Institution of Mechanical Engineers, Part J: Journal of Engineering Tribology},
	author = {Dobrica, Mihai and Fillon, Michel},
	month = feb,
	year = {2009},
	pages = {69--78}
}

@article{dobrica_reynolds_2005,
	title = {Reynolds' {Model} {Suitability} in {Simulating} {Rayleigh} {Step} {Bearing} {Thermohydrodynamic} {Problems}},
	volume = {48},
	doi = {10.1080/05698190500385088},
	language = {en},
	number = {4},
	journal = {Tribology Transactions},
	author = {Dobrica, Mihai B. and Fillon, Michel},
	month = aug,
	year = {2005},
	pages = {522--530},
}

@article{rahmani_analytical_2010,
	title = {An analytical approach for analysis and optimisation of slider bearings with infinite width parallel textures},
	volume = {43},
	doi = {10.1016/j.triboint.2010.02.016},
	number = {8},
	journal = {Tribology International},
	author = {Rahmani, R. and Mirzaee, I. and Shirvani, A. and Shirvani, H.},
	month = aug,
	year = {2010},
	pages = {1551--1565}
}

@article{armaly_experimental_1983,
	title = {Experimental and theoretical investigation of backward-facing step flow},
	volume = {127},
	doi = {10.1017/S0022112083002839},
	number = {1},
	journal = {Journal of Fluid Mechanics},
	author = {Armaly, B. F. and Durst, F. and Pereira, J. C. F. and Schönung, B.},
	month = feb,
	year = {1983},
	pages = {473},
}

@article{brown_applicability_1995,
	title = {Applicability of the {Reynolds} equation for modeling fluid flow between rough surfaces},
	volume = {22},
        doi = {10.1029/95GL02666},
	number = {18},
	journal = {Geophysical Research Letters},
	author = {Brown, R. Stephen and Stockman, Harlan W. and Reeves, Sally J.},
	month = sep,
	year = {1995},
	pages = {2537--2540},
}

@article{biswas_backward-facing_2004,
	title = {Backward-{Facing} {Step} {Flows} for {Various} {Expansion} {Ratios} at {Low} and {Moderate} {Reynolds} {Numbers}},
	volume = {126},
	doi = {10.1115/1.1760532},
	journal = {Journal of Fluids Engineering-Transactions of the Asme},
	author = {Biswas, G. and Breuer, M. and Durst, F.},
	month = may,
	year = {2004},
	pages = {362--374}
}

@article{gupta_new_2005,
	title = {A new paradigm for solving {Navier}–{Stokes} equations: streamfunction–velocity formulation},
	volume = {207},
	shorttitle = {A new paradigm for solving {Navier}–{Stokes} equations},
	doi = {10.1016/j.jcp.2005.01.002},
	number = {1},
	journal = {Journal of Computational Physics},
	author = {Gupta, M.M. and Kalita, J.C.},
	month = jul,
	year = {2005},
	pages = {52--68},
}

@misc{biswas_hoc_2017,
	title = {{HOC} simulation of {Moffatt} eddies and its flow topology in the triangular cavity flow},
	publisher = {arXiv: Computational Physics},
	author = {Biswas, S. and Kalita, J.C.},
        doi = {10.48550/arXiv.1710.06251},
	month = oct,
	year = {2017},
}

@article{marner_potential-velocity_2014,
	title = {On a potential-velocity formulation of {Navier}-{Stokes} equations},
	volume = {17},
	doi = {10.1134/S1029959914040110},
	journal = {Physical Mesomechanics},
	author = {Marner, Florian and Gaskell, Philip and Scholle, Markus},
	month = oct,
	year = {2014},
	pages = {124--130},
}

@article{sen_4oec_2015,
	title = {A {4OEC} scheme for the biharmonic steady {Navier}–{Stokes} equations in non-rectangular domains},
	volume = {196},
	doi = {10.1016/j.cpc.2015.05.024},
	journal = {Computer Physics Communications},
	author = {Sen, Shuvam and Kalita, Jiten C.},
	month = nov,
	year = {2015},
	pages = {113--133},
}

@article{takeuchi_extended_2019,
	title = {Extended {Reynolds} lubrication model for incompressible {Newtonian} fluid},
	volume = {4},
	issn = {2469-990X},
	doi = {10.1103/PhysRevFluids.4.114101},
	language = {en},
	number = {11},
	journal = {Physical Review Fluids},
	author = {Takeuchi, Shintaro and Gu, Jingchen},
	month = nov,
	year = {2019}
}

@book{leal_advanced_2007,
	title = {Advanced {Transport} {Phenomena}},
	language = {en},
	author = {Leal, L Gary},
	year = {2007},
	publisher = {Cambridge University Press},
}

@book{gallier_geometric_2011,
	address = {New York, NY},
	series = {Texts in {Applied} {Mathematics}},
	title = {Geometric {Methods} and {Applications}: {For} {Computer} {Science} and {Engineering}},
	volume = {38},
	isbn = {978-1-4419-9960-3 978-1-4419-9961-0},
	shorttitle = {Geometric {Methods} and {Applications}},
	language = {en},
	urldate = {2025-11-19},
	publisher = {Springer New York},
	author = {Gallier, Jean},
	year = {2011},
	doi = {10.1007/978-1-4419-9961-0},
}

@article{jain_numerically_2007,
	title = {Numerically {Stable} {Algorithms} for {Inversion} of {Block} {Tridiagonal} and {Banded} {Matrices}},
        url = {http://docs.lib.purdue.edu/ecetr/357},
	language = {en},
        journal = {ECE Technical Reports},
        number = {357},
        year = {2007},
	author = {Jain, Jitesh and Li, Hong and Cauley, Stephen and Koh, Cheng-Kok and Balakrishnan, Venkataramanan},
}

@article{he_modified_2023,
	title = {Modified {Finite} {Difference} {Methods} for {Reynold} {Equation} {With} {Film} {Thickness} {Discontinuity}},
	volume = {146},
	issn = {0742-4787},
	url = {https://doi.org/10.1115/1.4063442},
	doi = {10.1115/1.4063442},
	number = {024101},
	journal = {Journal of Tribology},
	author = {He, Qiang and Hu, Fengming and Huang, Weifeng and Hu, Yang and Cong, Guohui and Zhang, Yixun},
	month = oct,
	year = {2023},
}

@unpublished{dennis_separation_2026,
	   title = {Comparison of {Lubrication} {Theory} and {Stokes} {Flow} in {Corner} {Geometries} with {Flow} {Separation}},
      author={Dennis, Sarah and Fai, Thomas G.},
      year={2026},
      eprint={2501.18575},
      archivePrefix={arXiv},
      primaryClass={physics.flu-dyn},
      note = {arXiv:2501.18575 [physics]},
}
\bibliographystyle{siamplain}
\end{document}